\documentclass[12pt]{article}
\usepackage[T1]{fontenc}
\usepackage[utf8]{inputenc}
\usepackage[english]{babel}
\usepackage{amsmath, latexsym, amsthm, amsfonts, amssymb, amscd}
\usepackage{graphicx}
\usepackage{cite}
\usepackage{hyperref}
\usepackage{a4wide}
\usepackage{booktabs}
\usepackage{caption}
\graphicspath{{pictures_09_11_26/}}

\DeclareGraphicsExtensions{.pdf,.png,.jpg}

\newtheorem{theorem}{Theorem}
\newtheorem{lem}{Lemma}
\newtheorem{remark}{Remark}
\newtheorem{prop}{Proposition}

\begin{document}

\title{From oscillations to symmetry breaking in a two-dimensional
piecewise-smooth dynamical system}

\author{Oksana Satur\\[1mm]
{\small Institute of Mathematics of the National Academy of Sciences of Ukraine,}\\
{\small 3 Tereshchenkivska St., Kyiv, 01024, Ukraine}\\
{\small E-mail: \texttt{kseniajasko@gmail.com}}
}

\date{}

\maketitle

\begin{abstract}
We propose a discrete-time dynamical model in which two players, each distributed over two positions, interact only at the positions of their maximal presence, with an intensity set by an external resource. The main idea is that the evolution of the maxima closes into a self-contained two-dimensional piecewise-smooth map that can be analysed exactly. We show that a single critical level of the resource separates two outcomes: below it, the players exchange leadership forever along a stable cycle; above it, the symmetry between them breaks, and an arbitrarily small initial advantage grows into a permanent winner--loser hierarchy. The local stability of these regimes is proved, and their global attractivity is supported by numerical experiments. At the critical level itself, all dominance is erased and both players end in equal division of the resource. Mathematically, the transition is a border-collision bifurcation shaped by this symmetry.


\end{abstract}

\section{Introduction}

Many competitive processes are decided at the opponents' strongest positions rather than by their average state. In technology races a follower may overtake the leader (``leapfrogging''), so that leadership can alternate \cite{fudenberg1983}; in contests over several fronts each front is decided by the forces concentrated there \cite{roberson2006}; and in animal societies positive feedback can amplify small differences in competitive success into stable dominance hierarchies, with a transition between egalitarian and hierarchical states as the population density varies \cite{bonabeau1995}. The natural questions are whether leadership alternates or freezes, whether the competition can erase all advantages, and how these regimes are separated in parameter space; we answer them by a rigorous analysis of the attractors and of the transition between them, complemented by numerical evidence for the global attractivity of the asymmetric attractors.

A minimal mathematical setting for such leader-driven competition is a pair of probability distributions $\mathbf p^t$, $\mathbf r^t$ over two positions in which only the maximal components interact: in our model, at every step of discrete time $t$ the maximal component of $\mathbf p^t$ is multiplied by $\alpha(1-r^t_{\max})$ and that of $\mathbf r^t$ by $\alpha(1-p^t_{\max})$, where $\alpha>0$ is the external-resource parameter, the other components are left unchanged, and both vectors are renormalized. The pair of maxima $(x^t,y^t)=(p^t_{\max},r^t_{\max})$ then evolves autonomously by the \emph{map of maxima} $\bigl(x^{t+1},y^{t+1}\bigr)=\Phi\bigl(x^t,y^t\bigr)$, a continuous piecewise-smooth map on $\left(\tfrac12,1\right)^2$:
\begin{equation}\label{eq_macro_system}
x^{t+1} = \frac{\max\{\alpha x^t(1-y^t),\, 1-x^t\}}{\alpha x^t(1-y^t) + (1-x^t)}, \qquad
y^{t+1} = \frac{\max\{\alpha y^t(1-x^t),\, 1-y^t\}}{\alpha y^t(1-x^t) + (1-y^t)}.
\end{equation}
The non-smooth operator $\max$ determines whether the maximum of each distribution stays at the same position or switches to the other one.

The conceptual roots of modelling competitive interactions trace back to classical population dynamics, notably the continuous-time Lotka--Volterra equations and discrete models developed in the early 20th century \cite{hofbauer1998}. While these classical approaches primarily rely on scalar quantities to represent absolute population sizes, the approach initiated by Koshmanenko shifts the paradigm to a probabilistic framework: it merges the concepts of dynamical systems with competitive decision-making by describing the states of two or more opponents as stochastic vectors or general probability measures over a shared resource space. In this framework, the vector coordinates have a statistical meaning: they represent the probabilities of an opponent's presence within specific regions of a compartmentalized space. The temporal evolution of the system is governed by an interaction mapping, formalized via nonlinear difference equations, which redistributes these measures at every step of discrete time.

The foundation of this approach is a theorem proved by Koshmanenko, first for pairs of stochastic vectors \cite{kosh2003} and then for general probability measures \cite{kosh2004}: under the basic repulsive (alternative) interaction, every pair of non-identical, mutually non-singular distributions evolves towards a limiting pair of mutually singular measures, so that the opponents asymptotically divide the space into disjoint domains of exclusive presence. A systematic account of the spectral and structural aspects of this theory is given in the monograph \cite{koshbook}, and a recent survey of this approach in connection with opinion-formation models can be found in \cite{kosh2025}; the role of the point spectrum of equilibrium states in models of belief formation is discussed in \cite{koshsatur2025}. Over the past two decades this framework has generated a rich family of models~--- non-annihilating opponents \cite{albeverio2005}, cyclic migration of complex systems \cite{albeverio2008}, an abstract ``society'' with redistribution of social energy \cite{karataeva2020}, a mean-field ``game for power'' \cite{karataieva2019}, and the competition of alternative opponents in the presence of external support, which gives rise to compromise states \cite{karataieva2024}. Conceptually, these models belong to the same circle of problems as evolutionary game dynamics \cite{hofbauer1998}, the statistical physics of social dynamics \cite{castellano2009}, and bounded-confidence opinion formation \cite{hegselmann2002,marvel2012,shang2014}, but they are distinguished by an explicitly non-game-theoretic, measure-redistribution mechanism of interaction and by the possibility of a complete rigorous analysis.

Formally, the base model, a generalized discrete-time system of Lotka--Volterra type, describes the coordinate-wise evolution as follows:
\[
p_i^{t+1} = \frac{1}{z^t}p_i^t\big(1+\gamma\, r_i^t\big), \qquad r_i^{t+1} = \frac{1}{z^t}r_i^t\big(1+\gamma\, p_i^t\big),
\]
where $z^t=1+\gamma\,(\mathbf p^t,\mathbf r^t)$, with $(\mathbf p,\mathbf r)=\sum_i p_ir_i$, is the normalizing denominator ensuring stochasticity, common to both vectors, and $\gamma\neq0$ is the interaction parameter: $\gamma>0$ corresponds to attractive and $\gamma<0$ to repulsive interaction; in the latter case $\gamma\ge-1$, so that all numerators remain nonnegative, and $\gamma=-1$ gives the basic alternative interaction \cite{koshbook,albeverio2005}. Within this framework, two archetypes of interaction have been investigated in detail. Under \emph{repulsive} interaction ($\gamma<0$) the opponents suppress each other's presence at shared positions, which leads to spatial separation \cite{kosh2003,kosh2004,albeverio2005}. Under \emph{attractive} interaction ($\gamma>0$) the components of the distributions are mutually reinforced; here the main lines of research are the realization of a ``sure event''~--- convergence of a distribution to a vertex of the simplex~--- \cite{koshsatur2020,satur2021,koshkharchenko2017} and the study of the behaviour of trajectories, basins of attraction and fixed points \cite{satur2026}. Another direction was opened in \cite{satur2025}, where the interaction between two opponents is concentrated exclusively at the positions of their \emph{minimal} presence (``minimal players''): the external resource parameter stimulates or suppresses precisely the smallest components of the stochastic vectors, which corresponds to a risk-averse strategy of reinforcing one's weakest positions.

The present paper develops the opposite mechanism. For two positions, convergence to stationary states for $0<\alpha\le1$ and to the uniform distribution for $\alpha=2$ was proved in \cite{satur2025}, while the periodic regimes, a coexistence of a stationary and a $2$-periodic vector for $1<\alpha<2$ and a $4$-periodic regime for $\alpha>2$, were only observed numerically; here they occur in reversed order, with the same critical value $\alpha=2$ and the same collapse at it. Compared with \cite{satur2025}, the new results are the reduction to the map \eqref{eq_macro_system} with explicit switching curves; proofs of existence, uniqueness and local stability of the oscillatory and the symmetry-broken regimes for all $\alpha\neq2$; their closed-form values, involving the golden ratio and the plastic number; and the identification of the transition at $\alpha=2$ as a border collision accompanied by symmetry breaking. In applied terms the model describes the strategy of concentrating the main resources at the positions of maximal influence~--- direct competition of leaders rather than reinforcement of the weakest positions; such concentration mechanisms are typical for models of competition for resources and of resource exploitation in dynamic-game and discrete-time settings \cite{scheffran2000,cerboni2020}.

The map \eqref{eq_macro_system} isolates the pure deterministic mechanism of dominance switching. The model thereby joins two well-developed but rarely intersecting lines of research: the theory of max-type difference equations, where equations with the maximum operator are known to generate eventually periodic and highly structured dynamics \cite{groveladas2005,patula2004}, also for planar systems of max-type equations \cite{papaschinopoulos2003,stevic2013}, and the theory of piecewise-smooth maps and their border-collision bifurcations, originating from the works of Nusse and Yorke \cite{nusse1992,nusse1995} and Nusse, Ott and Yorke \cite{nusse1994}, extended to two-dimensional maps by Banerjee and Grebogi \cite{banerjee1999} (see also \cite{fournier2011}), and developed into a mature discipline in \cite{dibernardo2008,simpson2010,simpson2016,avrutin2019,sushko2025}. Organizing roles of switching manifolds in two-dimensional discontinuous maps \cite{avrutin2019,sushko2025}, the structure of stable and unstable sets of saddles in piecewise-smooth maps \cite{guney2026} and the bifurcation structure of chaotic attractors of two-dimensional piecewise-linear maps, including non-invertible ones \cite{ghosh2024}, are the subject of active current research; piecewise-smooth maps with border collisions also underlie applied models ranging from credit cycles \cite{sushko2016}, symmetric two-country models of trade and innovation \cite{sushko2018} and a three-region model of a developing economy \cite{commendatore2018} to renewable-resource exploitation \cite{cerboni2020}; bifurcation mechanisms of such two-dimensional nonsmooth maps are analysed, e.g., in \cite{panchuk2023}.

\section{Construction of the model}\label{model}

We consider two players $A$ and $B$ and a set of two positions $\Omega=\{\omega_1,\omega_2\}$. At $t=0$ the players are described by discrete probability distributions on $\Omega$,
\[
A\sim{\bf p}=(p_1,p_2),\qquad B\sim{\bf r}=(r_1,r_2),\qquad
p_i,\,r_i\ge0,\quad p_1+p_2=r_1+r_2=1,
\]
that is, by stochastic vectors whose coordinates $p_i$ and $r_i$ are the masses of $A$ and $B$, respectively, at the position $\omega_i$. A state of the system is a pair of stochastic vectors on $\Omega$; such a vector is called \emph{uniform} if it equals $\left(\tfrac12,\tfrac12\right)$, and a \emph{vertex} if it equals $(1,0)$ or $(0,1)$. At the moments $t=1,2,\dots$ of discrete time the players interact. We denote by $\divideontimes$ one step of this interaction, which transforms a state into the next one, so that the trajectory of the initial state is
\[
\{{\bf p},{\bf r}\}\equiv\{{\bf p}^0,{\bf r}^0\}
\stackrel{\divideontimes}{\longrightarrow}\{{\bf p}^1,{\bf r}^1\}
\stackrel{\divideontimes}{\longrightarrow}\cdots
\stackrel{\divideontimes}{\longrightarrow}\{{\bf p}^t,{\bf r}^t\}
\stackrel{\divideontimes}{\longrightarrow}\cdots,
\]
where ${\bf p}^t=(p_1^t,p_2^t)$ and ${\bf r}^t=(r_1^t,r_2^t)$ are stochastic vectors.

The interaction is concentrated at the positions where the masses of the players are maximal. If neither ${\bf p}^t$ nor ${\bf r}^t$ is uniform, the indices of these \emph{dominant positions},
\begin{equation*}
i^*=\arg\max_{i}p_i^t,\qquad j^*=\arg\max_{j}r_j^t,
\end{equation*}
are uniquely defined; we write $p^t_{\max}=p^t_{i^*}$ and $r^t_{\max}=r^t_{j^*}$ for the maximal values and call the other coordinates, $p_k^t$ with $k\neq i^*$ and $r_m^t$ with $m\neq j^*$, \emph{non-dominant}. The step $\divideontimes$ is defined by the evolution law
\begin{gather}\label{conf_max}
p^{t+1}_{i^*}=\frac{\alpha\, p^t_{\max}\,(1-r^t_{\max})}{z_p^t},\qquad
r^{t+1}_{j^*}=\frac{\alpha\, r^t_{\max}\,(1-p^t_{\max})}{z_r^t},\\
p_k^{t+1}=\frac{p_k^t}{z_p^t},\qquad
r_m^{t+1}=\frac{r_m^t}{z_r^t},\qquad
k\neq i^*,\quad m\neq j^*,\quad \alpha>0,\nonumber\\
z_p^t=\alpha\, p^t_{\max}\,(1-r^t_{\max})+p_k^t,\qquad
z_r^t=\alpha\, r^t_{\max}\,(1-p^t_{\max})+r_m^t,\nonumber
\end{gather}
for $t=0,1,\dots$; here $\alpha$ is the only parameter of the model, and the normalizing denominators $z_p^t$ and $z_r^t$ are the sums of the corresponding numerators, which ensures the stochasticity of ${\bf p}^{t+1}$ and ${\bf r}^{t+1}$. Thus, at every step the maximal coordinate of ${\bf p}^t$ is multiplied by $\alpha(1-r^t_{\max})$ (\emph{active interaction} at the dominant position), its non-dominant coordinate changes only through the normalization (\emph{passive redistribution}), and ${\bf r}^t$ is transformed in the same way with the factor $\alpha(1-p^t_{\max})$. In contrast with the base model of the Introduction, each vector has its own normalizing denominator, and the new value of each vector depends on the other vector only through the maximal value of the latter, irrespective of whether the dominant positions of the two vectors coincide. We call the dynamical system defined by the evolution law \eqref{conf_max} the \emph{full system}.

The right-hand sides of \eqref{conf_max} are uniquely defined at a state unless one of the vectors is uniform or both vectors are vertices. For a uniform vector the dominant position is not unique. Further, since $\alpha>0$ and $p^t_{\max}\ge\tfrac12$, the denominator $z_p^t$ vanishes only if $p_k^t=0$ and $r^t_{\max}=1$, i.e.\ only if both vectors are vertices; the same holds for $z_r^t$. A vertex is preserved by \eqref{conf_max}, since its non-dominant coordinate vanishes. If ${\bf p}^t$ is a vertex and ${\bf r}^t$ is not, then $1-p^t_{\max}=0$ and hence $r^{t+1}_{j^*}=0$, so that ${\bf r}^{t+1}$ is a vertex as well, and both normalizing denominators vanish at the next state; the case of a vertex ${\bf r}^t$ is symmetric. If neither vector is a vertex, then all numerators in \eqref{conf_max} are positive, and so are all coordinates of ${\bf p}^{t+1}$ and ${\bf r}^{t+1}$. Accordingly, we call a state \emph{degenerate} if at least one of its vectors is uniform or is a vertex, and \emph{non-degenerate} otherwise, and we consider only non-degenerate initial states $\{{\bf p}^0,{\bf r}^0\}$. Along the trajectory of such a state all coordinates remain positive, so that no vertex is ever reached, and the law \eqref{conf_max} determines the next state at every step at which the maximum of each vector is attained at a unique position.

The law \eqref{conf_max} refers to the positions only through the indices $i^*$ and $j^*$. Hence, if the coordinates of ${\bf p}^t$ are interchanged, then $i^*$ is replaced by the other index and $p^t_{\max}$ does not change, so that the coordinates of ${\bf p}^{t+1}$ are interchanged and ${\bf r}^{t+1}$ does not change; the same holds with the roles of ${\bf p}$ and ${\bf r}$ reversed. Moreover, the formulas for ${\bf p}^{t+1}$ and ${\bf r}^{t+1}$ are obtained from each other by interchanging ${\bf p}$ and ${\bf r}$, so that the interchange of ${\bf p}^t$ and ${\bf r}^t$ results in the interchange of ${\bf p}^{t+1}$ and ${\bf r}^{t+1}$. Thus each of these transformations commutes with the step $\divideontimes$, and a state preserved by such a transformation, or by a composition of them, is mapped to a state with the same property. Consequently, the set of states with ${\bf p}^t={\bf r}^t$ and the set of states with $p_1^t=r_2^t$ are invariant: the former consists of the states preserved by the interchange of the vectors, and the latter of the states preserved by this interchange combined with the interchange of the coordinates of both vectors.

We denote the maximal values by
\[
x^t:=p_{\max}^t,\qquad y^t:=r_{\max}^t,
\]
and call them the \emph{macroscopic parameters} of the system; for a non-degenerate state, $(x^t,y^t)\in\left(\tfrac12,1\right)^2$. A stochastic vector on $\Omega$ with a unique maximal coordinate is determined by the value and the position of its maximum, and, by the symmetries above, the positions of the maxima do not affect the new maximal values. Hence $x^{t+1}$ and $y^{t+1}$ are determined by $(x^t,y^t)$ alone. Explicitly, $x^{t+1}$ is the larger of the two coordinates of ${\bf p}^{t+1}$; since both coordinates have the same positive denominator $z_p^t$, the maximum can be taken over their numerators, that is, over the outcome $\alpha x^t(1-y^t)$ of the active interaction at the dominant position and the outcome $p_k^t=1-x^t$ of the passive redistribution at the other position. Similarly, $y^{t+1}$ is obtained by comparing $\alpha y^t(1-x^t)$ with $r_m^t=1-y^t$. In this way the pair of maxima evolves autonomously, by the map of maxima \eqref{eq_macro_system}.

In applied terms, the coordinates $p_i$ and $r_i$ are the probabilities of the presence of $A$ and $B$, respectively, at the position $\omega_i$, the dominant position of a player is the position of its maximal presence, and the parameter $\alpha$ is the level of an external resource. The law \eqref{conf_max} describes players that concentrate the interaction at the positions of their maximal presence: the maximal presence of each player is multiplied by a factor that is proportional to the external resource and to the presence of the other player outside the dominant position of the latter, whereas the presence at the other position changes only through the renormalization. A uniform vector describes a player without a dominant position, a vertex describes a player concentrated at a single position, and the macroscopic parameters measure the degree of dominance of the two players.

\section{The map of maxima and its attractors}\label{map}

For definiteness we write the evolution law for the case when the maxima of both players are attained at the first coordinate, $p_1^t>p_2^t$ and $r_1^t>r_2^t$ (so that $i^*=j^*=1$); all other arrangements of the dominant positions ($i^*=j^*=2$ or $i^*\neq j^*$) are obtained by relabelling the coordinates of each vector separately, and lead to the same formulas for the maxima.

In terms of the macroscopic parameters we have $p_1^t=x^t$ and $r_1^t=y^t$. By stochasticity ($p_1^t+p_2^t=1$, $r_1^t+r_2^t=1$), the second (non-dominant) coordinates are uniquely determined as $p_2^t=1-x^t$ and $r_2^t=1-y^t$.

Substituting these relations into the evolution law \eqref{conf_max}, we obtain the explicit system of difference equations:
\begin{gather}\label{conf2_max}
  p^{t+1}_{1} = \frac{\alpha x^t(1- y^t)}{z_p^t}, \qquad
  r^{t+1}_{1} = \frac{\alpha y^t(1- x^t)}{z_r^t}, \\
  p_2^{t+1} = \frac{1-x^t}{z_p^t}, \qquad
  r_2^{t+1} = \frac{1-y^t}{z_r^t}, \nonumber
\end{gather}
where, after substitution of the non-dominant coordinates, the normalizing denominators take the explicit algebraic form
$$z_p^t = \alpha x^t(1- y^t) + (1-x^t) = 1+ (\alpha-1) x^t -\alpha x^t y^t,$$
$$z_r^t = \alpha y^t(1- x^t) + (1-y^t) = 1+ (\alpha-1) y^t -\alpha x^t y^t.$$

Since each player has exactly one non-dominant coordinate, its numerator is exactly the complement of the current maximum to one ($1-x^t$ and $1-y^t$, respectively). This allows us to write the evolution of the macroscopic parameters explicitly as an autonomous two-dimensional system of nonlinear difference equations $\big(x^{t+1},y^{t+1}\big)=\Phi\big(x^t,y^t\big)$, which is precisely the map of maxima \eqref{eq_macro_system} displayed in the Introduction.
In what follows, for each coordinate of the map, the branch on which the maximum is attained at the first argument (e.g.\ $\alpha x(1-y)$) is called \emph{active}, since it corresponds to the dominance of the player's own mass; the branch on which the maximum is attained at the second argument (e.g.\ $1-x$) is called \emph{passive}.

In terms of the full system, an active coordinate means that the corresponding player keeps its dominant position, and a passive one that the position of its maximum switches to the opposite one. On the diagonal $x=y$ both coordinates are given by branches of the same type, passive for $\alpha x<1$ and active for $\alpha x>1$; we call these two cases the \emph{switching} and the \emph{preserving branch}, respectively. Since both arguments of each $\max$ are continuous and the denominators are positive on $\left(\tfrac12,1\right)^2$, the map $\Phi$ is continuous and piecewise smooth, with the switching curves
$$T_p=\{\alpha x(1-y)=1-x\},\qquad T_r=\{\alpha y(1-x)=1-y\},$$
on which the active and the passive numerators of the first, respectively second, coordinate coincide. On $T_p$ (respectively $T_r$) the first (respectively second) coordinate of $\Phi$ equals $\tfrac12$; thus $\Phi$ maps $\left(\tfrac12,1\right)^2\setminus(T_p\cup T_r)$ into $\left(\tfrac12,1\right)^2$, and the curves $T_p$, $T_r$ into the sides $x=\tfrac12$, $y=\tfrac12$ of the square. On each of the regions into which $T_p$ and $T_r$ divide the square, the coordinates of $\Phi$ are rational functions; across the switching curves $\Phi$ is continuous, but its first derivatives jump. This jump is a fold: since the passive value of each coordinate equals one minus the active one, the first coordinate of $\Phi$ equals $\tfrac12+\bigl|\tfrac{\alpha x(1-y)}{z_p}-\tfrac12\bigr|$, where $z_p=\alpha x(1-y)+1-x$, so that its gradient changes sign across $T_p$; similarly, the gradient of the second coordinate changes sign across $T_r$. In contrast, the evolution law \eqref{conf_max} has, in general, no continuous extension to states in which a vector has two equal maximal coordinates: the two possible choices of the dominant position give vectors that differ by the transposition of their coordinates. The maximal values of these vectors coincide, so that $\Phi$ does not inherit this discontinuity.

The presence of the non-smooth operator $\max$ in the numerators of \eqref{eq_macro_system} is the key mathematical feature of the proposed model. It is this operator that is responsible for the deterministic switching of the dominant positions during the evolution (in applied terms, a change of the ``leader'').

Since the macroscopic parameters belong to the open square $(x^t,y^t)\in\left(\tfrac12,1\right)^2$ (the boundaries are excluded: the value $\tfrac12$ corresponds to an undefined $\arg\max$, and the value $1$ to a degenerate vertex, from which both vectors become vertices in one step; the normalizing denominators vanish only under simultaneous concentration of both players), the structural behaviour of the system depends essentially on the positions at which the new maxima $x^{t+1}$ and $y^{t+1}$ are formed at the next step. In the two-dimensional state space two fundamental configurations naturally arise:
we say that a state $(\mathbf p^t,\mathbf r^t)$ is in \emph{Configuration~1} if
$i^*=j^*$, i.e.\ the maxima of both vectors are attained at the same
coordinate, and in \emph{Configuration~2} if $i^*\neq j^*$.

In applied terms, Configuration~1 corresponds to direct competition (the spatial priorities of the players coincide, and the interaction is applied
directly to the maximal concentrations), while Configuration~2 corresponds to avoidance of direct confrontation: the maximal influence of each player is directed to the position where the opponent's presence is minimal (cross interaction).

The specific nonlinearity of the evolution operator (in particular, the non-smooth function $\max$ in the
comparison of numerators in \eqref{eq_macro_system}) generates internal dynamics capable of deterministically switching the system between these two structural configurations.

We use the following names for the periodic orbits of the full system that arise below. For $x^*\in\left(\tfrac12,1\right)$, the orbit
$$\{(x^*, 1-x^*),\,(x^*, 1-x^*)\}\ \stackrel{\divideontimes}{\longleftrightarrow}\ \{(1-x^*, x^*),\,(1-x^*, x^*)\}$$
is called the \emph{synchronous $2$-cycle with value $x^*$}, and the orbit
$$\{(x^*, 1-x^*),\,(1-x^*, x^*)\}\ \stackrel{\divideontimes}{\longleftrightarrow}\
\{(1-x^*, x^*),\,(x^*, 1-x^*)\}$$
the \emph{mirror $2$-cycle with value $x^*$}. Along both orbits the position
of the maximum switches at every step, while its value is constant and equal
to $x^*$; the maxima of the two vectors are always at the same position in
the first case and at different positions in the second. For
$\tfrac12<y^*<x^*<1$, the orbit
$$\{(x^*,1-x^*),(1-y^*,y^*)\}\stackrel{\divideontimes}{\longrightarrow} \{(y^*,1-y^*),(x^*,1-x^*)\}\stackrel{\divideontimes}{\longrightarrow}$$
$$\stackrel{\divideontimes}{\longrightarrow} \{(1-x^*,x^*),(y^*,1-y^*)\}\stackrel{\divideontimes}{\longrightarrow}
\{(1-y^*,y^*),(1-x^*,x^*)\}\stackrel{\divideontimes}{\longrightarrow}\cdots$$
is called the \emph{$4$-cycle with values $(x^*,y^*)$}. Along it the maxima of
the two vectors take the values $x^*$ and $y^*$ and exchange them at every
step, while the state alternates between Configurations~2 and~1.

Denote by $\sigma(a,b)=(b,a)$ the permutation of coordinates. Formulas
\eqref{eq_macro_system} are symmetric under the interchange $x\leftrightarrow
y$, so $\Phi$ is equivariant \cite{golubitsky2002}: $\Phi\circ\sigma=\sigma\circ\Phi$.
Consequently, if a point $w=(x,y)\in\left(\tfrac12,1\right)^2$ with $x\neq y$ satisfies
$\Phi(w)=\sigma w$, then $\Phi(\sigma w)=\sigma\Phi(w)=w$, so that
$\{w,\sigma w\}$ is a $2$-cycle of $\Phi$ with the spatio-temporal symmetry
$\sigma$ in the terminology of \cite{golubitsky2002}: applying $\sigma$ to the
cycle is the same as advancing it by one step. We call such a cycle an
\emph{exchange $2$-cycle}; along it the maxima of the two vectors swap their
values at every step. The exchange $2$-cycle $\{(x^*,y^*),(y^*,x^*)\}$ of
$\Phi$ with $y^*<x^*$ and the $4$-cycle of the full system with values
$(x^*,y^*)$ are the same object, written in the macroscopic parameters and in
the full coordinates, respectively; general $2$-cycles of $\Phi$, with
$\Phi^2(w)=w$ and $\Phi(w)\neq\sigma w$, are not exchange $2$-cycles.

Note that initial distributions of the form $(1,0)$ and $(0,1)$ are excluded from the phase space: their non-dominant coordinate vanishes, so the corresponding vector stays at the vertex, while
the other vector is sent to a vertex in one step; at the resulting state of simultaneous
concentration of both players the normalizing denominators vanish and the law
\eqref{conf2_max} loses its meaning. In the space of macroscopic parameters a vertex
corresponds to a side $x=1$ or $y=1$ of the closed square, which $\Phi$ maps to the corner
$(1,1)$; this corner is never reached in finite time from the interior of the square, since both numerators in \eqref{eq_macro_system} are positive in the interior, so that $x^{t+1},y^{t+1}<1$.

In contrast, the diagonally opposite corner $\left(\tfrac12,\tfrac12\right)$ of the boundary (equal division of the resource) plays a central role in what follows.

\begin{figure}[htbp]
\centering
\includegraphics[width=0.77\textwidth]{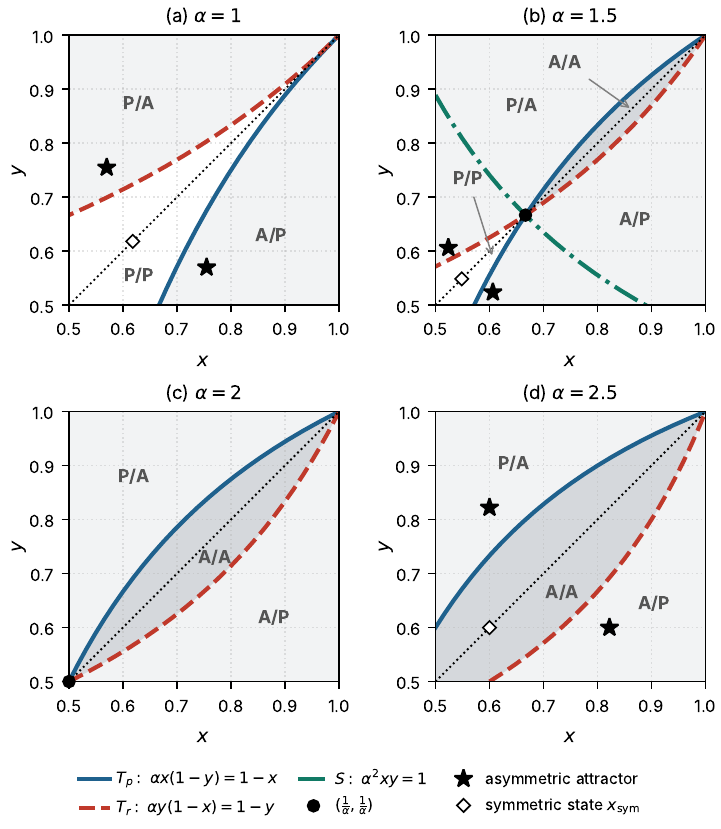}
\caption{Switching curves $T_p$ and $T_r$, the curve $S$ and the diagonal
(dotted) in the phase square $\left(\tfrac12,1\right)^2$. The labels A/A, A/P, P/A, P/P
indicate whether the first and the second coordinate of $\Phi$ are given by
the active (A) or the passive (P) branch (shading: A/A dark, A/P and P/A light,
P/P white); the curve $S$ meets the square only for $1<\alpha<2$. The dot is
the point $\left(\tfrac1\alpha,\tfrac1\alpha\right)$, where $T_p$ and $T_r$ cross; the
stars are the exchange $2$-cycle (Theorem~\ref{thm_cycle_general}) for
$\alpha<2$ and the asymmetric fixed points (Theorem~\ref{thm_fixed_general})
for $\alpha>2$, and the diamond is the symmetric state
$(x_{\mathrm{sym}},x_{\mathrm{sym}})$. (a) $\alpha=1$: $T_p$ and $T_r$ meet
only at the corner $(1,1)$, and the region A/A is empty, since double activity
would require $\alpha^2xy>1$; the coordinates of the diamond and of the stars
are given by the golden ratio and the plastic number
(Remark~\ref{rem_constants}). (b) $\alpha=1.5$: the three curves meet at the
point $\left(\tfrac1\alpha,\tfrac1\alpha\right)$ of the diagonal. (c) $\alpha=2$ and (d)
$\alpha=2.5$: the region P/P is empty, since double passivity would require
$\alpha^2xy<1$, whereas $\alpha^2xy>\tfrac{\alpha^2}{4}\ge1$ in the whole
square; at $\alpha=2$ the point $\left(\tfrac1\alpha,\tfrac1\alpha\right)$ is the corner
$\left(\tfrac12,\tfrac12\right)$, and for $\alpha>2$ the curves $T_p$ and $T_r$ do not
meet inside the square.}
\label{fig_phase_space_comparison}
\end{figure}

\begin{remark}\label{rem_exceptional}
The law \eqref{conf_max} becomes undefined if at some step one of the vectors becomes uniform, $\left(\tfrac12,\tfrac12\right)$: then the $\arg\max$ is ambiguous. This happens exactly when the active and passive numerators of the corresponding player coincide at the previous step, i.e.\ when the previous state lies on one of the switching curves $T_p$, $T_r$ (see Fig.~\ref{fig_phase_space_comparison}).
 In addition, for $\alpha\in(1,2)$ the curve $$S=\{\alpha^2xy=1\}$$ without its point
$\left(\tfrac1\alpha,\tfrac1\alpha\right)\in T_p\cap T_r$ is mapped onto the diagonal $\{x=y\}$
in one step. Indeed, the ratios $\tfrac{\alpha x(1-y)}{1-x}$ and $\tfrac{\alpha y(1-x)}{1-y}$
of the active to the passive numerators have the product $\alpha^2xy$ and coincide only for
$x=y$, so that at the other points of $S$ one coordinate of the map is active and the other
passive (or vice versa); for such a combination of regimes the equality of the new maxima $\tfrac{\alpha x(1-y)}{z_p}=\tfrac{1-y}{z_r}$ is equivalent to $\alpha x z_r=z_p$, which after expansion gives $\alpha^2xy(1-x)=1-x$, i.e.\ $\alpha^2xy=1$ (the symmetric combination of regimes gives the same curve). Once on the diagonal, the trajectory of $\Phi$ remains on it forever by the equivariance of $\Phi$, and therefore it cannot converge to an attractor lying off the diagonal. Note that the curve $S$ is empty both for $\alpha\le1$, since $\alpha^2xy\le xy<1$, and for $\alpha\ge2$, since $\alpha^2xy>\tfrac{\alpha^2}{4}\ge1$.
\end{remark}

Denote
\begin{equation}\label{eq_exceptional}
E_\alpha=\bigcup_{t\ge0}\Phi^{-t}\bigl(T_p\cup T_r\cup S\bigr).
\end{equation}
This is at most a countable union of smooth curves. Indeed, on each of the four combinations of regimes (active or passive numerator for each coordinate) the map $\Phi$ is smooth, and a direct computation gives for the Jacobian of each branch
$$\det J_\Phi(x,y)=\pm\,\frac{\alpha^2(1-x)(1-y)(1-xy)}{z_p^2\,z_r^2}\neq0, \qquad(x,y)\in\left(\tfrac12,1\right)^2$$
(the sign is positive when both coordinates are in regimes of the same type and negative for mixed regimes), so each branch is a local diffeomorphism and the preimage of a smooth curve under $\Phi$ is again an at most countable union of smooth arcs. Since preimages of Lebesgue null sets under local diffeomorphisms are null sets, $E_\alpha$ has Lebesgue measure zero and, in particular, no interior points.
All statements about the asymmetric dynamics below are formulated for initial states whose pair of maxima $(x^0,y^0)$ lies outside $E_\alpha$; in particular, on the diagonal the curve $T_p$ reduces to the single point $x=y=\tfrac1\alpha$ (nonempty only for $1<\alpha<2$); since moreover $S=\emptyset$ for $\alpha\notin(1,2)$, the set $E_\alpha$ does not meet the diagonal for $\alpha\notin(1,2)$, and the symmetric statements are automatically covered by the same condition.

We first record a property of the branch structure of $\Phi$ that will be
used below. A step of a trajectory of $\Phi$ is called \emph{doubly active}
if both coordinates of \eqref{eq_macro_system} are given by the active
branches.

\begin{lem}\label{lem_double_active}
Let $0<\alpha<2$, and let a trajectory of $\Phi$ lie in
$\left(\tfrac12,1\right)^2$ outside the exceptional set $E_\alpha$ of
\eqref{eq_exceptional}. Then
\begin{enumerate}
\item for $0<\alpha\le1$ the trajectory has no doubly active steps;
\item for $1<\alpha<2$ every series of consecutive doubly active steps of the
trajectory is finite.
\end{enumerate}
\end{lem}

\begin{proof}
\emph{The case $0<\alpha\le1$.} Double activity means that both
inequalities
$$\alpha x(1-y)>1-x,\qquad \alpha y(1-x)>1-y$$
hold. Both sides of each inequality are positive, so they can be multiplied
and divided by $(1-x)(1-y)>0$; this gives $\alpha^2xy>1$, which is
impossible, since $\alpha^2xy\le xy<1$.

\emph{The case $1<\alpha<2$: monotonicity along a series.} Suppose, to the
contrary, that all steps from some $t_0$ on are doubly active. Then
$x^{t+1}=\dfrac{\alpha x^t(1-y^t)}{z_p^t}$ for all $t\ge t_0$, and
$$x^{t+1}-x^t=\frac{x^t\bigl[\alpha(1-y^t)-z_p^t\bigr]}{z_p^t}
=\frac{x^t(1-x^t)\bigl[\alpha(1-y^t)-1\bigr]}{z_p^t},$$
where we used $\alpha(1-y^t)-z_p^t=\alpha(1-y^t)(1-x^t)-(1-x^t)$. Since
$y^t>\tfrac12$, we have $\alpha(1-y^t)<\tfrac\alpha2<1$, so $x^{t+1}<x^t$;
symmetrically, $x^t>\tfrac12$ implies $y^{t+1}<y^t$.

\emph{Passage to the limit.} Both sequences are decreasing and bounded
below by $\tfrac12$, hence they have limits
$\bar x,\bar y\in\bigl[\tfrac12,1\bigr)$; the limits are strictly less than
one, since $x^t\le x^{t_0}<1$ and $y^t\le y^{t_0}<1$. All functions in the
recurrence $x^{t+1}z_p^t=\alpha x^t(1-y^t)$ are continuous, so in the limit
$\bar x\,\bar z_p=\alpha\bar x(1-\bar y)$, and since $\bar x\ge\tfrac12>0$,
we get $\bar z_p=\alpha(1-\bar y)$, i.e.
$$\alpha\bar x(1-\bar y)+(1-\bar x)=\alpha(1-\bar y)
\iff(1-\bar x)\bigl[1-\alpha(1-\bar y)\bigr]=0.$$
Since $\bar x<1$, this yields $\alpha(1-\bar y)=1$, i.e.
$\bar y=1-\tfrac1\alpha<\tfrac12$ for $\alpha<2$, a contradiction with
$\bar y\ge\tfrac12$.
\end{proof}

The next lemma describes the action of $\Phi$ in the logarithmic
coordinates $\ell(s)=\ln\frac{s}{1-s}$, $s\in[\tfrac12,1)$, of the maxima.

\begin{lem}\label{lem_order}
Let $\alpha>0$, $(x,y)\in[\tfrac12,1)^2$, let $(x',y')$ be given by
\eqref{eq_macro_system}, and let $\ell(s)=\ln\frac{s}{1-s}$. Then
\begin{enumerate}
\item $|\ell(x')-\ell(y')|\le2|\ln\alpha|+2\ln2$ and
$\ell(x')+\ell(y')\le\max\bigl\{2|\ln\alpha|+2\ln2,\ 2|\ell(x)-\ell(y)|\bigr\}$;
\item if $\alpha\ge2$ and $x>y$, then $\alpha x(1-y)>1-x$ and $x'>y'$;
\item if $\alpha>2$, then
$|\ell(x')-\ell(y')|\ge\min\bigl\{2\ln\tfrac\alpha2,\ \tfrac32\,|\ell(x)-\ell(y)|\bigr\}$.
\end{enumerate}
\end{lem}

\begin{proof}
Dividing the numerators and the denominators in \eqref{eq_macro_system} by
$1-x$ and $1-y$, respectively, we get
$$x'=\chi(A),\qquad y'=\chi(B),\qquad
A=\frac{\alpha x(1-y)}{1-x},\quad B=\frac{\alpha y(1-x)}{1-y},\qquad
\chi(t)=\frac{\max\{t,1\}}{1+t}.$$
Since $\ell(\chi(t))=|\ln t|$ for $t>0$, we have
$\ell(x')=|\ln A|$ and $\ell(y')=|\ln B|$. Moreover,
$$AB=\alpha^2xy,\qquad \frac AB=\frac{x\bigl(1-y\bigr)^2}{y\bigl(1-x\bigr)^2}=e^{m(x)-m(y)},
\qquad m(s)=\ln\frac{s}{\bigl(1-s\bigr)^2}.$$
For real $a$, $b$ one has $|a|+|b|=\max\{|a+b|,|a-b|\}$ and
$\bigl||a|-|b|\bigr|=\min\{|a+b|,|a-b|\}$. Taking $a=\ln A$ and $b=\ln B$, we
obtain
\begin{equation}\label{eq_log_maxima}
\ell(x')+\ell(y')=\max\{|P|,|M|\},\qquad
|\ell(x')-\ell(y')|=\min\{|P|,|M|\},
\end{equation}
where $P=\ln\bigl(\alpha^2xy\bigr)$ and $M=m(x)-m(y)$. Since
$\frac{dm}{ds}=\frac{1+s}{s(1-s)}=(1+s)\frac{d\ell}{ds}$ and
$1+s\in[\tfrac32,2)$, the mean value theorem applied to $m\circ\ell^{-1}$
gives
$$\tfrac32\,|\ell(x)-\ell(y)|\le|M|\le2\,|\ell(x)-\ell(y)|.$$
Further, $\ln x+\ln y\in[-2\ln2,0)$, so that $|P|\le2|\ln\alpha|+2\ln2$. By
\eqref{eq_log_maxima}, $|\ell(x')-\ell(y')|=\min\{|P|,|M|\}\le2|\ln\alpha|+2\ln2$
and $\ell(x')+\ell(y')=\max\{|P|,|M|\}\le\max\bigl\{2|\ln\alpha|+2\ln2,\ 2|\ell(x)-\ell(y)|\bigr\}$.

Let $\alpha\ge2$ and $x>y$. Then $AB=\alpha^2xy>\tfrac{\alpha^2}{4}\ge1$ and
$A/B>1$, since $x>y$ and $1-y>1-x$. Hence $A>B$, $A>1/B$ and $A>1$; the last
inequality means $\alpha x(1-y)>1-x$. The function $\chi$ strictly increases
on $[1,\infty)$ and $\chi(t)=\chi(1/t)$, so that
$y'=\chi(B)=\chi(\max\{B,1/B\})<\chi(A)=x'$.

If $\alpha>2$, then $|P|=\ln\bigl(\alpha^2xy\bigr)\ge\ln\tfrac{\alpha^2}{4}=2\ln\tfrac\alpha2$,
and by \eqref{eq_log_maxima} and the lower bound for $|M|$ we obtain
$|\ell(x')-\ell(y')|=\min\{|P|,|M|\}\ge\min\bigl\{2\ln\tfrac\alpha2,\ \tfrac32|\ell(x)-\ell(y)|\bigr\}$.
\end{proof}

Lemma~\ref{lem_order} describes one step of $\Phi$; along trajectories it
has the following consequences. Let $c_\alpha=2|\ln\alpha|+2\ln2$. By the
first inequality of Lemma~\ref{lem_order}(1),
$|\ell(x^{t-1})-\ell(y^{t-1})|\le c_\alpha$ for $t\ge2$, and then the second
inequality of Lemma~\ref{lem_order}(1) gives
$\ell(x^t)+\ell(y^t)\le\max\{c_\alpha,2c_\alpha\}=4|\ln\alpha|+4\ln2$ for
every $t\ge2$. Thus every trajectory enters after two steps the compact set
$$K_\alpha=\bigl\{(x,y)\in[\tfrac12,1)^2:\ \ell(x)+\ell(y)\le4|\ln\alpha|+4\ln2\bigr\}$$
and remains there. By Lemma~\ref{lem_order}(2) and induction on $t$, for
$\alpha\ge2$ the inequality $x^t>y^t$, once valid at $t=0$, holds for all
$t$, and the maximum of the corresponding vector keeps its position at every
step. Finally, let $\alpha>2$, $x^0\neq y^0$ and
$d_t=|\ell(x^t)-\ell(y^t)|$, so that $d_0>0$. By Lemma~\ref{lem_order}(3),
$d_{t+1}\ge\tfrac32d_t$ as long as $\tfrac32d_t<2\ln\tfrac\alpha2$; hence
$d_t\ge2\ln\tfrac\alpha2$ for some $t$, and by the same inequality this
bound is preserved at all subsequent steps.

In the logarithmic coordinates of Lemma~\ref{lem_order}, the map $\Phi$ is a
fold of a smooth injective map whose Jacobian matrix does not depend on
$\alpha$.

\begin{lem}\label{lem_logform}
Let $\alpha>0$, let $\ell(s)=\ln\frac{s}{1-s}$, and for $u,v\ge0$ let
\begin{gather*}
\eta(u,v)=\ln\alpha+u-\ln\bigl(1+e^{v}\bigr),\\
H(u,v)=\bigl(\eta(u,v),\,\eta(v,u)\bigr),\qquad
\Psi(u,v)=\bigl(|\eta(u,v)|,\,|\eta(v,u)|\bigr).
\end{gather*}
Then
\begin{enumerate}
\item for $(x,y)\in[\tfrac12,1)^2$ and $(x',y')=\Phi(x,y)$ one has
$\bigl(\ell(x'),\ell(y')\bigr)=\Psi\bigl(\ell(x),\ell(y)\bigr)$;
\item the first coordinate of $\Phi(x,y)$ is given by the active branch if
$\eta(\ell(x),\ell(y))>0$ and by the passive branch if
$\eta(\ell(x),\ell(y))<0$; the same holds for the second coordinate and
$\eta(\ell(y),\ell(x))$;
\item the map $H$ is smooth and injective, and
$$J_H(u,v)=\begin{pmatrix}1&-y\\-x&1\end{pmatrix},\qquad
x=\ell^{-1}(u),\quad y=\ell^{-1}(v);$$
\item on each set on which $\tau_1=\operatorname{sign}\eta(u,v)$ and
$\tau_2=\operatorname{sign}\eta(v,u)$ are constant and nonzero,
$$\Psi=\begin{pmatrix}\tau_1&0\\0&\tau_2\end{pmatrix}H,$$
so that $\Psi$ is injective on this set and
$|\det J_\Psi(u,v)|=1-xy\le\tfrac34$.
\end{enumerate}
\end{lem}

\begin{proof}
Put $u=\ell(x)$ and $v=\ell(y)$; then $\frac{x}{1-x}=e^{u}$ and
$1-y=\frac{1}{1+e^{v}}$. For the quantities
$A=\frac{\alpha x(1-y)}{1-x}$ and $B=\frac{\alpha y(1-x)}{1-y}$ from the proof
of Lemma~\ref{lem_order} this gives
$\ln A=\ln\alpha+u-\ln\bigl(1+e^{v}\bigr)=\eta(u,v)$ and, by symmetry,
$\ln B=\eta(v,u)$. As shown in that proof, $\ell(x')=|\ln A|$ and
$\ell(y')=|\ln B|$, which is the formula for $\Psi$. The inequality $A>1$ is
equivalent to $\alpha x(1-y)>1-x$, i.e.\ the first coordinate is given by
the active branch, and $A<1$ means that it is given by the passive one; the
same holds for $B$ and the second coordinate.

Let $\kappa(t)=\ln(1+e^{t})$. Then $\kappa'(t)=\frac{e^{t}}{1+e^{t}}\in(0,1)$,
$\kappa'(u)=x$ and $\kappa'(v)=y$, which gives $J_H$. If
$H(u_1,v_1)=H(u_2,v_2)$, then $u_1-u_2=\kappa(v_1)-\kappa(v_2)$ and
$v_1-v_2=\kappa(u_1)-\kappa(u_2)$. By the mean value theorem,
$|\kappa(a)-\kappa(b)|<|a-b|$ for $a\neq b$; hence $v_1\neq v_2$ would give
$|u_1-u_2|<|v_1-v_2|\le|u_1-u_2|$, which is impossible. Thus $v_1=v_2$, then
$u_1=u_2$, and $H$ is injective. On a set where $\tau_1$ and $\tau_2$ are
constant and nonzero, $|\eta(u,v)|=\tau_1\eta(u,v)$ and
$|\eta(v,u)|=\tau_2\eta(v,u)$, which gives the representation of $\Psi$, its
injectivity on this set and $|\det J_\Psi|=|\det J_H|=1-xy$. Finally,
$x,y\ge\tfrac12$ implies $1-xy\le\tfrac34$.
\end{proof}

\subsection{The oscillatory range $0<\alpha<2$}

We consider first symmetric and mirror-symmetric initial states and then asymmetric ones.

\begin{prop}\label{prop_sym_general}
Let $0<\alpha<2$, let $x_{\mathrm{sym}}=\frac{-1+\sqrt{1+4\alpha}}{2\alpha}$ be
the unique root of $\alpha x^2+x-1=0$ in $\left(\tfrac12,1\right)$, and let the initial
state not belong to the exceptional set $E_\alpha$ of \eqref{eq_exceptional}.
Then
\begin{enumerate}
\item if $\mathbf p^0=\mathbf r^0$, the system converges to the synchronous
$2$-cycle with value $x_{\mathrm{sym}}$;
\item if $p_1^0=r_2^0$, the system converges to the mirror $2$-cycle with
value $x_{\mathrm{sym}}$.
\end{enumerate}
\end{prop}

\begin{proof}
By the symmetry of the formulas for the macroscopic parameters, the diagonal $x^t=y^t$ is invariant: for $x^t=y^t=x$ both new coordinates are equal to $\max\{\alpha x(1-x),\,1-x\}/z$, where $z=(1-x)(\alpha x+1)$. Cancelling $(1-x)>0$, we obtain the one-dimensional map $$x^{t+1}=g(x^t)=\frac{\max\{\alpha x^t,\,1\}}{\alpha x^t+1}.$$
The map $g$ sends $\left(\tfrac12,1\right)\setminus\left\{\tfrac1\alpha\right\}$ into $\left(\tfrac12,1\right)$: on the switching branch $\alpha x<1$ we have $g(x)=\tfrac{1}{\alpha x+1}>\tfrac12$, on the preserving branch $\alpha x>1$ we have $g(x)=\tfrac{\alpha x}{\alpha x+1}>\tfrac12$, and in both cases $g(x)<1$, while $g\left(\tfrac1\alpha\right)=\tfrac12$; the kink point $x=\tfrac1\alpha$ (relevant only for $\alpha\in(1,2)$) and its preimages belong to the exceptional set $E_\alpha$ of \eqref{eq_exceptional}. On each branch
$$|g'(x)|=\frac{\alpha}{\bigl(\alpha x+1\bigr)^2}\le\frac{\alpha}{\left(\tfrac\alpha2+1\right)^2} =\frac{4\alpha}{\bigl(\alpha+2\bigr)^2}=:k<1,$$
since $\bigl(\alpha+2\bigr)^2-4\alpha=\alpha^2+4>0$. Since $g$ is continuous and piecewise smooth, it is Lipschitz with constant $k$ on the whole of $\left(\tfrac12,1\right)$. The fixed point lies on the switching branch: the equation $x=\tfrac{1}{\alpha x+1}$ gives $\alpha x^2+x-1=0$, whence $x_{\mathrm{sym}}=\tfrac{-1+\sqrt{1+4\alpha}}{2\alpha}$, and the branch condition $\alpha x_{\mathrm{sym}}<1$ is equivalent to $\alpha<2$. Then for every admissible $x_0\in\left(\tfrac12,1\right)$
$$\bigl|x^{t}-x_{\mathrm{sym}}\bigr|=\bigl|g(x^{t-1})-g(x_{\mathrm{sym}})\bigr|
\le k\,\bigl|x^{t-1}-x_{\mathrm{sym}}\bigr|\le\cdots\le k^{t}\bigl|x_0-x_{\mathrm{sym}}\bigr|\to0,$$ i.e.\ the iterations converge to $x_{\mathrm{sym}}$ geometrically. Since $\alpha x_{\mathrm{sym}}<1$, the position of the maximum switches at every step, so that the system converges to the synchronous $2$-cycle with value $x_{\mathrm{sym}}$ for $\mathbf p^0=\mathbf r^0$ and to the mirror one for $p_1^0=r_2^0$.
\end{proof}

We now turn to asymmetric initial states, for which the symmetric $2$-cycle is replaced by a $4$-cycle.

\begin{theorem}\label{thm_cycle_general}
Let $0<\alpha<2$, let $x^\ast$ be the unique root of the equation
\begin{equation}\label{eq_cubic_general}
\alpha^{3}x^{3}+\alpha x^{2}+(1-\alpha)x-1=0
\end{equation}
in the interval $\left(\tfrac12,1\right)$, and let $y^\ast=\tfrac{1}{1+\alpha x^\ast}$,
so that $\tfrac12<y^\ast<x^\ast<1$. Then
\begin{enumerate}
\item the pair $(x^\ast,y^\ast)$, $(y^\ast,x^\ast)$ is an exchange $2$-cycle of
$\Phi$; in the full
system \eqref{conf_max} it corresponds, up to relabelling of positions, to the $4$-cycle with values $(x^\ast,y^\ast)$;
\item this cycle is locally asymptotically stable; its multipliers are
complex conjugate (a stable focus) with common modulus
$$1-x^\ast y^\ast=\frac{1+(\alpha-1)x^\ast}{1+\alpha x^\ast}<\tfrac34,$$
so every trajectory entering a sufficiently small neighbourhood of
$(x^\ast,y^\ast)$ converges to the cycle geometrically with asymptotic
ratio $1-x^\ast y^\ast$ per two steps.
\end{enumerate}
\end{theorem}

\begin{proof}
 We look for a solution of
$\Phi(x,y)=(y,x)$ with $\tfrac12<y<x<1$ in which the first coordinate of
$\Phi$ is active and the second one is passive. The exchange
condition then reads
\begin{align}
\frac{1}{z_p}\,\alpha x(1-y)=y,&\qquad z_p=\alpha x(1-y)+(1-x),
\label{eq_exch_p}\\
\frac{1}{z_r}\,(1-y)=x,&\qquad z_r=\alpha y(1-x)+(1-y).
\label{eq_exch_r}
\end{align}
Multiplying \eqref{eq_exch_r} by $z_r$, we obtain
$1-y=\alpha xy(1-x)+x(1-y)$, i.e.\ $(1-x)(1-y)=\alpha xy(1-x)$, whence
$1-y=\alpha xy$ and
\begin{equation}\label{eq_exch_y}
y=\frac{1}{1+\alpha x}.
\end{equation}
Multiplying \eqref{eq_exch_p} by $z_p$, we obtain
$\alpha x(1-y)=\alpha xy(1-y)+y(1-x)$, whence
\begin{equation}\label{eq_exch_sq}
\alpha x\bigl(1-y\bigr)^2=y(1-x).
\end{equation}
Since $1-y=\alpha xy$, we have $\bigl(1-y\bigr)^2=\alpha^2x^2y^2$, and
\eqref{eq_exch_sq} becomes $\alpha^3x^3y^2=y(1-x)$, i.e.\
$y=\dfrac{1-x}{\alpha^3x^3}$. Equating this with \eqref{eq_exch_y}, we get
$$\frac{1}{1+\alpha x}=\frac{1-x}{\alpha^3x^3}
\iff \alpha^3x^3=1+\alpha x-x-\alpha x^2,$$
which is equation \eqref{eq_cubic_general}.

Let $$P(x)=\alpha^3x^3+\alpha x^2+(1-\alpha)x-1.$$ Then $P(1)=\alpha^3>0$ and $$P\left(\tfrac12\right)=\tfrac{\alpha^3}{8}-\tfrac{\alpha}{4}-\tfrac12=:h_0(\alpha).$$ We have $h_0(0)=-\tfrac12$, $h_0(2)=0$ and $h_0'(\alpha)=\tfrac{3\alpha^2}{8}-\tfrac14$, so $h_0$ decreases on $\big(0,\sqrt{2/3}\big)$ and increases on $\big(\sqrt{2/3},2\big)$; hence $h_0(\alpha)<0$ on $(0,2)$, and $P\left(\tfrac12\right)<0<P(1)$. Moreover, $P'(x)=3\alpha^3x^2+2\alpha x+(1-\alpha)$ is a quadratic with vertex at $-\tfrac{1}{3\alpha^2}<0$, so for $x\ge\tfrac12$ $$P'(x)\ge P'\left(\tfrac12\right)=\tfrac{3\alpha^3}{4}+1>0.$$ Therefore $P$ is strictly increasing on $\left[\tfrac12,1\right]$ and has a unique root $x^\ast\in\left(\tfrac12,1\right)$.

We show that $\alpha x^\ast<1$. For $0<\alpha\le1$ this is obvious. For $1<\alpha<2$ we have $\tfrac1\alpha\in\left(\tfrac12,1\right)$ and $P\left(\tfrac1\alpha\right)=\tfrac2\alpha-1>0=P(x^\ast)$, whence $x^\ast<\tfrac1\alpha$. Consequently $1<1+\alpha x^\ast<2$, and \eqref{eq_exch_y} gives $\tfrac12<y^\ast<1$.

The inequality $y^\ast<x^\ast$ is equivalent to $x^\ast(1+\alpha x^\ast)>1$, i.e.\ $\alpha \bigl(x^\ast\bigr)^2+x^\ast-1>0$. The equation $\alpha x^2+x-1=0$ has the positive root $x_{\mathrm{sym}}=\tfrac{-1+\sqrt{1+4\alpha}}{2\alpha}\in\left(\tfrac12,1\right)$ for $0<\alpha<2$, and $\alpha x^2+x-1>0$ is equivalent to $x>x_{\mathrm{sym}}$. Using $\alpha x_{\mathrm{sym}}^2=1-x_{\mathrm{sym}}$, hence $\alpha^3x_{\mathrm{sym}}^3=\alpha^2x_{\mathrm{sym}}(1-x_{\mathrm{sym}})$ and $\alpha(1-x_{\mathrm{sym}})=\alpha^2x_{\mathrm{sym}}^2$, we find
\begin{align*}
P(x_{\mathrm{sym}})&=\alpha^2x_{\mathrm{sym}}(1-x_{\mathrm{sym}})
+\alpha x_{\mathrm{sym}}^2+(1-\alpha)x_{\mathrm{sym}}-1
=\alpha^2x_{\mathrm{sym}}(1-x_{\mathrm{sym}})-\alpha x_{\mathrm{sym}}\\
&=\alpha x_{\mathrm{sym}}\bigl(\alpha(1-x_{\mathrm{sym}})-1\bigr)
=\alpha x_{\mathrm{sym}}(\alpha x_{\mathrm{sym}}-1)(\alpha x_{\mathrm{sym}}+1).
\end{align*}
Since $\alpha x_{\mathrm{sym}}=\tfrac{-1+\sqrt{1+4\alpha}}{2}<1$ for $\alpha<2$, we get $P(x_{\mathrm{sym}})<0=P(x^\ast)$, so $x^\ast>x_{\mathrm{sym}}$ and $y^\ast<x^\ast$.

It remains to verify the choice of branches. We use two identities of the cycle: \eqref{eq_exch_y} gives $1-y^\ast=\alpha x^\ast y^\ast$, and \eqref{eq_cubic_general}, rewritten as $\alpha^3\big(x^\ast\big)^3=(1-x^\ast)(1+\alpha x^\ast)$, gives $1-x^\ast=\alpha^3\big(x^\ast\big)^3 y^\ast$. The first coordinate is active:
$$\alpha x^\ast(1-y^\ast)=\alpha^2\bigl(x^\ast\bigr)^2y^\ast \;>\;\alpha^3\bigl(x^\ast\bigr)^3y^\ast=1-x^\ast,$$
since the inequality is equivalent to $\alpha x^\ast<1$. The second coordinate is passive:
$$1-y^\ast=\alpha x^\ast y^\ast\;>\;\alpha^4\bigl(x^\ast\bigr)^3\bigl(y^\ast\bigr)^2 =\alpha y^\ast(1-x^\ast),$$
since after cancelling $\alpha x^\ast y^\ast>0$ the inequality becomes $\alpha^3\bigl(x^\ast\bigr)^2y^\ast<1$, and by the second identity $\alpha^3\bigl(x^\ast\bigr)^2y^\ast=\tfrac{1-x^\ast}{x^\ast}<1$, because $x^\ast>\tfrac12$.

Thus $\{(x^\ast,y^\ast),(y^\ast,x^\ast)\}$ is an exchange $2$-cycle of
$\Phi$. Along it the maxima of the full system pass successively through
the positions $(1,2)$, $(1,1)$, $(2,1)$, $(2,2)$ and return to $(1,2)$ after
four steps, which is the $4$-cycle with values $(x^\ast,y^\ast)$.

 Since the regime inequalities at the points
of the orbit are strict, $\Phi$ is smooth in neighbourhoods of these points.
The stability of the $2$-cycle $(x^\ast,y^\ast)\leftrightarrow(y^\ast,x^\ast)$
is determined by the monodromy matrix, the Jacobian of $\Phi$ over
the full period:
$$M=J_\Phi(y^\ast,x^\ast)\,J_\Phi(x^\ast,y^\ast).$$
The cycle is locally asymptotically stable if both eigenvalues of $M$ (the
multipliers of the cycle) lie in the open unit disc
\cite[Theorem~4.33]{elaydi}. In what follows we drop the asterisks, write
$x:=x^\ast$, $y:=y^\ast$, and use the identities of the cycle obtained
above:
\begin{equation}\label{eq_cycle_id}
1-y=\alpha xy,\qquad 1-x=\alpha^3x^3y,\qquad 1+\alpha x=\tfrac1y.
\end{equation}
At the cycle $x'=y$ and $y'=x$, i.e.\ $\tfrac{\alpha x(1-y)}{z_p}=y$ and
$\tfrac{1-y}{z_r}=x$, whence
$$z_p=\frac{\alpha x(1-y)}{y}=\alpha^2x^2,\qquad
z_r=\frac{1-y}{x}=\alpha y .$$

In a neighbourhood of the cycle the first coordinate is given by the active
branch $x'=\alpha x(1-y)/z_p$ and the second by the passive branch
$y'=(1-y)/z_r$. By the quotient rule, the numerators of the partial
derivatives simplify identically: for $\partial x'/\partial x$ we have
$z_p-x\bigl(\alpha(1-y)-1\bigr)\equiv1$, so the numerator equals
$\alpha(1-y)$; for $\partial x'/\partial y$ it equals
$-\alpha x\bigl(z_p-\alpha x(1-y)\bigr)=-\alpha x(1-x)$; for the second
coordinate, the numerator $1-y$ does not depend on $x$, and
$\partial z_r/\partial y=\alpha(1-x)-1$. Substituting
$z_p=\alpha^2x^2$, $z_r=\alpha y$ and \eqref{eq_cycle_id}, we obtain
\begin{align*}
\frac{\partial x'}{\partial x}&=\frac{\alpha(1-y)}{z_p^2}
=\frac{1-y}{\alpha^3x^4}=\frac{y}{\alpha^2x^3}, &
\frac{\partial x'}{\partial y}&=-\frac{\alpha x(1-x)}{z_p^2}
=-\frac{1-x}{\alpha^3x^3}=-y,\\[1mm]
\frac{\partial y'}{\partial x}&=\frac{\alpha y(1-y)}{z_r^2}
=\frac{1-y}{\alpha y}=x, &
\frac{\partial y'}{\partial y}&=\frac{-z_r-(1-y)\bigl(\alpha(1-x)-1\bigr)}{z_r^2}.
\end{align*}
In the last expression the numerator equals
$$-\alpha y-\alpha xy\bigl(\alpha(1-x)-1\bigr)
=-\alpha y\bigl((1-x)+\alpha x(1-x)\bigr)=-\alpha y(1-x)(1+\alpha x),$$
so, applying \eqref{eq_cycle_id} twice,
$$\frac{\partial y'}{\partial y}=-\frac{(1-x)(1+\alpha x)}{\alpha y}
=-\frac{1-x}{\alpha y^2}=-\frac{\alpha^2x^3}{y}.$$
Hence the one-step Jacobian matrix is
$$A:=J_\Phi(x^\ast,y^\ast)=
\begin{pmatrix} a & -y \\ x & b\end{pmatrix},\qquad
a=\frac{y}{\alpha^2x^3},\quad b=-\frac{\alpha^2x^3}{y},\quad ab=-1 .$$
By the equivariance $\Phi\circ\sigma=\sigma\circ\Phi$, differentiation gives
$J_\Phi(y^\ast,x^\ast)=SAS$, where
$S=\bigl(\begin{smallmatrix}0&1\\1&0\end{smallmatrix}\bigr)$; multiplication
by $S$ on both sides interchanges the rows and the columns:
$$B:=J_\Phi(y^\ast,x^\ast)=\begin{pmatrix} b & x \\ -y & a\end{pmatrix}.$$
Multiplying and using $ab=-1$, we get
$$M=BA=\begin{pmatrix} x^2-1 & b(x-y) \\ a(x-y) & y^2-1\end{pmatrix},$$
whence
\begin{equation}\label{eq_trM_general}
\operatorname{Tr}M=x^2+y^2-2,
\end{equation}
and, since $\det A=\det B=ab+xy=-(1-xy)$,
\begin{equation}\label{eq_detM_general}
\det M=\bigl(1-xy\bigr)^2.
\end{equation}

 By \cite[(4.3.8)]{elaydi}, both multipliers lie in the
open unit disc if and only if
$$1-\det M>0,\qquad 1+\operatorname{Tr}M+\det M>0,\qquad
1-\operatorname{Tr}M+\det M>0.$$
By \eqref{eq_trM_general}--\eqref{eq_detM_general},
\begin{align*}
1-\det M&=xy\,(2-xy),\\
1+\operatorname{Tr}M+\det M&=\bigl(x-y\bigr)^2+\bigl(xy\bigr)^2,\\
1-\operatorname{Tr}M+\det M&=4-\bigl(x+y\bigr)^2+\bigl(xy\bigr)^2,
\end{align*}
and all three expressions are strictly positive for $x,y\in\left(\tfrac12,1\right)$:
$0<xy<1$, $\bigl(xy\bigr)^2>0$, and $\bigl(x+y\bigr)^2<4$. Hence, for every $\alpha\in(0,2)$,
both multipliers lie in the open unit disc, and the cycle is locally
asymptotically stable.

Moreover, the discriminant of the characteristic equation factorizes as
\begin{align*}
D&=\bigl(\operatorname{Tr}M\bigr)^2-4\det M
=\bigl(\operatorname{Tr}M-2(1-xy)\bigr)\bigl(\operatorname{Tr}M+2(1-xy)\bigr)\\
&=\bigl(\bigl(x+y\bigr)^2-4\bigr)\bigl(x-y\bigr)^2=-(2-x-y)(2+x+y)\bigl(x-y\bigr)^2<0,
\end{align*}
since $x+y<2$ and $x\neq y$. Therefore the multipliers are complex
conjugate (a stable focus) on the whole range $0<\alpha<2$, and their
common modulus is
$$\rho(M)=\sqrt{\det M}=1-x^\ast y^\ast
=\frac{1+(\alpha-1)x^\ast}{1+\alpha x^\ast}<\tfrac34,$$
where the last inequality follows from $x^\ast y^\ast>\tfrac14$; this is a
uniform in $\alpha$ bound for the multipliers. By the theorem on stability
in the first approximation \cite[Theorem~4.33]{elaydi}, every trajectory
from a sufficiently small neighbourhood of $(x^\ast,y^\ast)$ converges to
the cycle; since the multipliers are complex conjugate, in a suitable basis
$M$ is $1-x^\ast y^\ast$ times a rotation matrix, so that the deviation is
multiplied asymptotically by $1-x^\ast y^\ast$ per two steps.
\end{proof}

\begin{remark}\label{rem_constants}
For $\alpha=1$ Proposition~\ref{prop_sym_general} and
Theorem~\ref{thm_cycle_general} take a particularly simple form. The
symmetric value is the root of $x^2+x-1=0$, i.e.\
$x_{\mathrm{sym}}=\tfrac{\sqrt5-1}{2}=\varphi^{-1}$, where
$\varphi=\tfrac{1+\sqrt5}{2}$ is the golden ratio; on the diagonal
$g(x)=\frac1{1+x}$, and the iterates of $g$ starting from $x=1$ are the ratios
$F_k/F_{k+1}$ of consecutive Fibonacci numbers. The cubic
\eqref{eq_cubic_general} becomes $x^3+x^2-1=0$, whose root is
$x^\ast=\psi^{-1}\approx0.7549$, where $\psi\approx1.3247$, the real root of
$\psi^3-\psi-1=0$, is the plastic number; moreover,
$y^\ast=\frac1{1+x^\ast}=\bigl(x^\ast\bigr)^2=\psi^{-2}\approx0.5698$, and the modulus
of the multipliers of the cycle equals $1-x^\ast y^\ast=y^\ast$. Thus the
passage from the $2$-cycle to the $4$-cycle is accompanied by the passage
from the quadratic constant associated with the Fibonacci sequence to its
cubic analogue associated with the Padovan sequence
\cite{stewart1996,aarts2001}.
\end{remark}

\begin{lem}\label{lem_list_2cycles}
Let $\alpha>0$ and $x_{\mathrm{sym}}(\alpha)=\tfrac{-1+\sqrt{1+4\alpha}}{2\alpha}$. Then
\begin{enumerate}
\item for $0<\alpha<2$ the solutions of the equation $\Phi(x,y)=(y,x)$ in
$\left(\tfrac12,1\right)^2$ are the diagonal point $x=y=x_{\mathrm{sym}}(\alpha)$ and the
points $(x^\ast,y^\ast)$, $(y^\ast,x^\ast)$ of Theorem~\ref{thm_cycle_general};
\item for $\alpha>2$ the only solution of this equation in $\left(\tfrac12,1\right)^2$ is
the diagonal point $x=y=1-\tfrac1\alpha$;
\item for $\alpha=2$ this equation has no solutions in $\left(\tfrac12,1\right)^2$.
\end{enumerate}
\end{lem}

\begin{proof}
Each coordinate of $\Phi$ is given by one of two numerators (active or
passive), so it suffices to examine the four combinations of regimes.

\emph{Both coordinates active.} Then
$\tfrac{\alpha x(1-y)}{z_p}=y$ and $\tfrac{\alpha y(1-x)}{z_r}=x$.
Multiplying the first equation by $z_p$ and regrouping (as in the derivation
of~\eqref{eq_exch_sq}), we get $\alpha x\bigl(1-y\bigr)^2=y(1-x)$; symmetrically,
$\alpha y\bigl(1-x\bigr)^2=x(1-y)$. Dividing the first relation by the second,
$$x^2\bigl(1-y\bigr)^3=y^2\bigl(1-x\bigr)^3\iff h(x)=h(y),\qquad h(u):=\frac{u^2}{\bigl(1-u\bigr)^3}.$$
The function $h$ is strictly increasing on $(0,1)$, so $x=y$, and then
$\alpha(1-x)=1$, i.e.\ $x=1-\tfrac1\alpha$. This point lies in
$\left(\tfrac12,1\right)$ if and only if $\alpha>2$, and there the activity condition
$\alpha x>1$ holds.

\emph{Both coordinates passive.} Then $\tfrac{1-x}{z_p}=y$ and
$\tfrac{1-y}{z_r}=x$. The first equation gives
$1-x=y\bigl[\alpha x(1-y)+(1-x)\bigr]$, i.e.\ $(1-x)(1-y)=\alpha xy(1-y)$,
whence $1-x=\alpha xy$; symmetrically $1-y=\alpha xy$. Hence $x=y$ and
$\alpha x^2+x-1=0$, i.e.\ $x=x_{\mathrm{sym}}(\alpha)$. The inequality
$x_{\mathrm{sym}}>\tfrac12$ is equivalent to $\sqrt{1+4\alpha}>1+\alpha$,
i.e.\ to $\alpha<2$.

\emph{First coordinate active, second passive.} This is system~\eqref{eq_exch_p}--\eqref{eq_exch_r}. Its reduction in the proof of
Theorem~\ref{thm_cycle_general} to $y=\tfrac{1}{1+\alpha x}$ and $P(x)=0$ is
purely algebraic and holds for every $\alpha>0$. Moreover,
$P'(x)\ge P'\left(\tfrac12\right)=\tfrac{3\alpha^3}{4}+1>0$ for $x\ge\tfrac12$, so $P$
is strictly increasing on $[\tfrac12,\infty)$ for every $\alpha>0$. For
$0<\alpha<2$ the unique solution is $(x^\ast,y^\ast)$, as shown in that
proof. For $\alpha\ge2$ we have
$P\left(\tfrac12\right)=h_0(\alpha)\ge0$, since $h_0(2)=0$ and
$h_0'(\alpha)=\tfrac{3\alpha^2}{8}-\tfrac14>0$ for $\alpha>\sqrt{2/3}$;
hence $P(x)>0$ on $\left(\tfrac12,1\right)$, and there are no solutions of this type.

\emph{First coordinate passive, second active.} The substitution
$(x,y)\mapsto(y,x)$ together with the equivariance
$\Phi\circ\sigma=\sigma\circ\Phi$ reduces this case to the previous one: for
$0<\alpha<2$ the solution is $(y^\ast,x^\ast)$, the second point of the same
cycle, and for $\alpha\ge2$ there are no solutions.
\end{proof}

In particular, the exchange $2$-cycle of Theorem~\ref{thm_cycle_general} is
unique up to permutation of coordinates.

\subsection{Case $\alpha=2$: global collapse}

\begin{lem}\label{lem_Gmono}
Let $\alpha=2$, $(x,y)\in[\tfrac12,1)^2$ and $(x',y')=\Phi(x,y)$. Then
\begin{enumerate}
\item $\max(x',y')<\max(x,y)$ if $\min(x,y)>\tfrac12$;
\item $\max(x',y')=\max(x,y)$ and $\min(x',y')>\tfrac12$ if
$\min(x,y)=\tfrac12<\max(x,y)$.
\end{enumerate}
\end{lem}

\begin{proof}
Formulas \eqref{eq_macro_system} define the values of the new maxima also on
the boundary of the square: if two coordinates of a vector are equal, the
values of the new maxima do not depend on which of them is taken as
dominant. Since the statement of the lemma and the map $\Phi$ are symmetric
with respect to the permutation $\sigma$, we may assume without loss of
generality that $x\ge y$, so that $\tfrac12\le y\le x<1$ and $x>\tfrac12$.
By \eqref{eq_macro_system} with $\alpha=2$, the vectors at the next step are
$${\mathbf p}' = \left(\frac{2x(1-y)}{z_p},\ \frac{1-x}{z_p}\right),\qquad
{\mathbf r}' = \left(\frac{2y(1-x)}{z_r},\ \frac{1-y}{z_r} \right),$$
where $z_p=1+x-2xy>0$ and $z_r=1+y-2xy>0$. Since $\max(x',y')$ is the
largest of these four quantities, it suffices to compare each of them
with $x$.

\emph{First coordinate of $\mathbf p'$.} We have $\dfrac{2x(1-y)}{z_p}\le x$, with equality only for $y=\tfrac12$.
Indeed, multiplying by $z_p$ and dividing by $x>0$, we get $2(1-y)\le1+x-2xy$, and
$$1+x-2xy-2(1-y)=(x-1)+2y(1-x)=(1-x)(2y-1)\ge0,$$
since $1-x>0$ and $2y-1\ge0$.

\emph{Second coordinate of $\mathbf p'$.} We have $\dfrac{1-x}{z_p}<x$. Indeed, multiplying by $z_p$, we get
$1-x<x(1+x-2xy)$, i.e.\ $2x+x^2-2x^2y-1>0$. Since $y\le x$, we have
$-2x^2y\ge-2x^3$, hence
$$2x+x^2-2x^2y-1\ \ge\ 2x+x^2-2x^3-1=(2x-1)(1-x)(x+1)>0 .$$

\emph{First coordinate of $\mathbf r'$.} We have $\dfrac{2y(1-x)}{z_r}<y$. Indeed, multiplying by $z_r$ and dividing by
$y>0$, we get $2(1-x)<1+y-2xy$, and
$$1+y-2xy-2(1-x)=(y-1)+2x(1-y)=(1-y)(2x-1)>0 .$$
Hence $\dfrac{2y(1-x)}{z_r}<y\le x$.

\emph{Second coordinate of $\mathbf r'$.} We have $\dfrac{1-y}{z_r}<x$. Indeed, multiplying by $z_r$, we get
$1-y<x(1+y-2xy)$, i.e.\ $x+xy-2x^2y-1+y>0$. Using
$1+x-2x^2=(1-x)(1+2x)$, we obtain
$$(x-1)+y\bigl(1+x-2x^2\bigr)=(1-x)\bigl[y(1+2x)-1\bigr]>0,$$
since $y(1+2x)\ge\tfrac12(1+2x)>1$.

If $y>\tfrac12$, then all four coordinates are strictly less
than $x$, whence $\max(x',y')<\max(x,y)$.

Let $y=\tfrac12$. Then the first coordinate of $\mathbf p'$ equals $x$, the other three are less than $x$, and therefore
$\max(x',y')=x=\max(x,y)$. Moreover, $z_p=1$, so the maximum of
$\mathbf p'$ equals $x>\tfrac12$, while $z_r=\tfrac32-x<1$, so the maximum
of $\mathbf r'$ equals $\dfrac{1/2}{\tfrac32-x}>\tfrac12$. Hence
$\min(x',y')>\tfrac12$.
\end{proof}

If at some step one of the vectors is uniform, the dominant position in
\eqref{conf_max} is chosen arbitrarily (for instance, the previous one is
kept); this choice does not affect the values of the maxima.

\begin{theorem}\label{thm_alpha2}
Let $\alpha=2$. Then for every pair of non-degenerate initial stochastic
vectors $\mathbf p^0,\mathbf r^0$ the trajectory of the full system with evolution \eqref{conf_max} converges to the uniform distributions:
$$\lim_{t\to\infty}\mathbf p^t=\left(\frac12,\frac12\right),\qquad
\lim_{t\to\infty}\mathbf r^t=\left(\frac12,\frac12\right).$$
\end{theorem}

\begin{proof}
 The maximum of a two-component
stochastic vector is not less than $\tfrac12$, so
$w^t:=(x^t,y^t)\in[\tfrac12,1)^2$ and $G_t:=\max(x^t,y^t)\ge\tfrac12$ for all $t$. If
$w^t=\left(\tfrac12,\tfrac12\right)$ for some $t$, then $z_p=z_r=1$ and the trajectory
remains at this point, so the statement holds. Otherwise, by
Lemma~\ref{lem_Gmono}, $G_{t+1}\le G_t$ for all $t$, and the sequence $G_t$
has a limit $G^\ast\ge\tfrac12$.

Suppose that
$G^\ast>\tfrac12$. All points $w^t$ lie in the compact square
$[\tfrac12,G_0]^2$, so there is a subsequence $w^{t_k}\to\bar w=(\bar x,\bar y)$,
and $\max(\bar x,\bar y)=G^\ast>\tfrac12$. The map $\Phi$ is continuous on
$[\tfrac12,1)^2$, since both arguments of each $\max$ in
\eqref{eq_macro_system} are continuous and the denominators are positive.
Denote $(\bar x',\bar y')=\Phi(\bar w)$ and $(\bar x'',\bar y'')=\Phi^2(\bar w)$.
Then $w^{t_k+2}=\Phi^2\bigl(w^{t_k}\bigr)\to\Phi^2(\bar w)$, and
$$\max(\bar x'',\bar y'')=\lim_{k\to\infty}G_{t_k+2}=G^\ast=\max(\bar x,\bar y).$$
On the other hand, apply Lemma~\ref{lem_Gmono} to the trajectory starting
at $\bar w$. If $\min(\bar x,\bar y)>\tfrac12$, then
$\max(\bar x',\bar y')<\max(\bar x,\bar y)$, and
$\max(\bar x'',\bar y'')\le\max(\bar x',\bar y')$. If
$\min(\bar x,\bar y)=\tfrac12$, then
$\max(\bar x',\bar y')=\max(\bar x,\bar y)$ and
$\min(\bar x',\bar y')>\tfrac12$, so
$\max(\bar x'',\bar y'')<\max(\bar x',\bar y')$. In both cases
$\max(\bar x'',\bar y'')<\max(\bar x,\bar y)$, a contradiction. Hence
$G^\ast=\tfrac12$.

Thus $\tfrac12\le\min(x^t,y^t)\le\max(x^t,y^t)=G_t\to\tfrac12$, whence
$x^t\to\tfrac12$ and $y^t\to\tfrac12$. Since the coordinates of
$\mathbf p^t$ are $x^t$ and $1-x^t$, and those of $\mathbf r^t$ are $y^t$
and $1-y^t$, we obtain $\mathbf p^t\to\left(\tfrac12,\tfrac12\right)$ and
$\mathbf r^t\to\left(\tfrac12,\tfrac12\right)$.
\end{proof}

\begin{remark}\label{rem_alpha2_cycle}
By Lemma~\ref{lem_list_2cycles}(3), for $\alpha=2$ the map $\Phi$ has no
exchange $2$-cycles. Moreover, since $P'(x)\ge1$ for $x\ge\tfrac12$ and $P(x^\ast)=0$, the
values of the $4$-cycle of Theorem~\ref{thm_cycle_general} satisfy
$0<x^\ast(\alpha)-\tfrac12\le-P\left(\tfrac12\right)=-h_0(\alpha)\to0$ as
$\alpha\to2^-$, and hence $y^\ast(\alpha)=\tfrac{1}{1+\alpha x^\ast(\alpha)}\to\tfrac12$. More precisely, regarding $P$ as a function of $(x,\alpha)$, at the point
$\left(\tfrac12,2\right)$ we have $P=0$, $\partial P/\partial x=7$ and
$\partial P/\partial\alpha=\tfrac54$; by the implicit function theorem,
$$x^\ast(\alpha)=\tfrac12+\tfrac5{28}(2-\alpha)+O\bigl(\bigl(2-\alpha\bigr)^2\bigr),\qquad
y^\ast(\alpha)=\tfrac12+\tfrac1{28}(2-\alpha)+O\bigl(\bigl(2-\alpha\bigr)^2\bigr)$$
as $\alpha\to2^-$, so that $x^\ast-y^\ast=\tfrac17(2-\alpha)+O\bigl(\bigl(2-\alpha\bigr)^2\bigr)$.
Thus, as $\alpha\to2^-$, the $4$-cycle shrinks linearly in $2-\alpha$ to the
corner $\left(\tfrac12,\tfrac12\right)$, which at $\alpha=2$ becomes the
global attractor.
\end{remark}

In applied terms, Theorem~\ref{thm_alpha2} describes the collapse of the
competition: at the critical level of the external resource neither
player retains dominance, and the system arrives at the state of complete
uncertainty, the equal division of the resource between the positions.

\subsection{Case $\alpha>2$: symmetry breaking and the winner--loser fixed point}

For $\alpha>2$ the external resource is strong enough for a dominant
coordinate to be reinforced: by \eqref{eq_macro_system}, an active
coordinate grows whenever $\alpha(1-y)>1$, which is now compatible with
$y>\tfrac12$. As a consequence, the maxima may keep their positions, and
the oscillatory regimes of the range $0<\alpha<2$ give way to fixed points.

\begin{prop}\label{prop_sym_fixed}
Let $\alpha>2$. Then
\begin{enumerate}
\item if $\mathbf p^0=\mathbf r^0$, the system converges to the fixed point
$$\Bigl\{\left(1-\tfrac1\alpha,\tfrac1\alpha\right),\ \left(1-\tfrac1\alpha,\tfrac1\alpha\right)\Bigr\};$$
\item if $p_1^0=r_2^0$, the system converges to the fixed point
$$\Bigl\{\left(1-\tfrac1\alpha,\tfrac1\alpha\right),\ \left(\tfrac1\alpha,1-\tfrac1\alpha\right)\Bigr\}.$$
\end{enumerate}
\end{prop}

\begin{proof}
As in the proof of Proposition~\ref{prop_sym_general}, the equality of the values of the maxima is preserved at every step, and their common value evolves by the map
$$x^{t+1}=g(x^t),\qquad g(x)=\frac{\max\{\alpha x,\,1\}}{\alpha x+1}.$$
Since $\alpha>2$, we have $\alpha x>\tfrac\alpha2>1$ for every $x\in\left(\tfrac12,1\right)$, so that $g(x)=\tfrac{\alpha x}{\alpha x+1}$ on the whole interval; in particular, the maximum of each vector keeps its position at every step. The map $g$ sends $\left(\tfrac12,1\right)$ into itself, and
$$0<g'(x)=\frac{\alpha}{\bigl(\alpha x+1\bigr)^2}<\frac{\alpha}{\left(\tfrac\alpha2+1\right)^2}
=\frac{4\alpha}{\bigl(\alpha+2\bigr)^2}=:k<1 .$$
The fixed point equation $x=\tfrac{\alpha x}{\alpha x+1}$ gives $x=1-\tfrac1\alpha\in\left(\tfrac12,1\right)$. By the mean value theorem, $\bigl|x^t-\left(1-\tfrac1\alpha\right)\bigr|\le k^t\bigl|x^0-\left(1-\tfrac1\alpha\right)\bigr|\to0$. Since the positions of the maxima do not change, for $\mathbf p^0=\mathbf r^0$ both maxima remain at the same position, and for $p_1^0=r_2^0$ at different positions, which gives the two fixed points indicated.
\end{proof}

We now turn to asymmetric initial states.

\begin{theorem}\label{thm_fixed_general}
Let $\alpha>2$, let
\begin{equation}\label{eq_fixed}
x^\ast=1-\frac{1}{\alpha\bigl(\alpha-1\bigr)^2},\qquad y^\ast=1-\frac{1}{\alpha},
\end{equation}
and let
\begin{equation}\label{eq_rho}
\rho(\alpha)=\sqrt{\frac{\alpha^2-\alpha+1}{\alpha^2(\alpha-1)}} .
\end{equation}
Then
\begin{enumerate}
\item $\tfrac12<y^\ast<x^\ast<1$, and the points $(x^\ast,y^\ast)$ and
$(y^\ast,x^\ast)$ are the only fixed points of $\Phi$ in $\left(\tfrac12,1\right)^2$
off the diagonal;
\item these fixed points are locally asymptotically stable; the Jacobian
$J$ of $\Phi$ at $(x^\ast,y^\ast)$ has zero trace and real eigenvalues
$\pm\rho(\alpha)$, and $J^2=\rho^2(\alpha)I$; the function $\rho$ strictly
decreases on $(2,\infty)$, with $\rho\bigl(2^+\bigr)=\tfrac{\sqrt3}{2}$ and
$\rho(\alpha)\to0$ as $\alpha\to\infty$; every trajectory entering a
sufficiently small neighbourhood of $(x^\ast,y^\ast)$ converges to it
geometrically with asymptotic ratio $\rho^2(\alpha)$ per two steps;
\item along such a trajectory of the full system \eqref{conf_max}, the
position of the maximum of the vector whose maximal value tends to
$x^\ast$ is eventually constant, and this vector converges to
$(x^\ast,1-x^\ast)$ up to relabelling of positions, whereas the position of
the maximum of the other vector changes at every step, and this vector
converges to the $2$-periodic regime
$(y^\ast,1-y^\ast)\leftrightarrow(1-y^\ast,y^\ast)$.
\end{enumerate}
\end{theorem}

\begin{proof}
We look for a fixed point of $\Phi$ in which the first coordinate is active
and the second one is passive:
\begin{align}
\frac{\alpha x(1-y)}{z_p}=x,&\qquad z_p=\alpha x(1-y)+(1-x),
\label{eq_fix_p}\\
\frac{1-y}{z_r}=y,&\qquad z_r=\alpha y(1-x)+(1-y).
\label{eq_fix_r}
\end{align}
Dividing \eqref{eq_fix_p} by $x>0$, we get $z_p=\alpha(1-y)$, i.e.\
$\alpha x(1-y)+(1-x)=\alpha(1-y)$, whence $(1-x)\bigl(1-\alpha(1-y)\bigr)=0$.
Since $x<1$, this gives $\alpha(1-y)=1$, i.e.\ $y=y^\ast$. Multiplying
\eqref{eq_fix_r} by $z_r$, we get
$1-y=\alpha y^2(1-x)+y(1-y)$, i.e.\ $\bigl(1-y\bigr)^2=\alpha y^2(1-x)$, whence
$$1-x=\frac{\bigl(1-y\bigr)^2}{\alpha y^2}=\frac{1/\alpha^2}{\alpha\bigl(\alpha-1\bigr)^2/\alpha^2}
=\frac{1}{\alpha\bigl(\alpha-1\bigr)^2},$$
i.e.\ $x=x^\ast$. Thus system \eqref{eq_fix_p}--\eqref{eq_fix_r} has the
unique solution $(x^\ast,y^\ast)$.

This point is admissible: $y^\ast>\tfrac12$ since $\alpha>2$;
$x^\ast>\tfrac12$ since $\alpha\bigl(\alpha-1\bigr)^2$ is increasing and equals $2$ at
$\alpha=2$; and $x^\ast>y^\ast$ since $\alpha\bigl(\alpha-1\bigr)^2>\alpha$ for
$\alpha>2$. At the fixed point $z_p=\alpha(1-y^\ast)=1$ and
$z_r=\tfrac{1-y^\ast}{y^\ast}=\tfrac{1}{\alpha-1}$. The first coordinate is
active, since $\alpha x^\ast(1-y^\ast)=x^\ast>1-x^\ast$. The second
coordinate is passive: by $\bigl(1-y^\ast\bigr)^2=\alpha \bigl(y^\ast\bigr)^2(1-x^\ast)$, the
inequality $1-y^\ast>\alpha y^\ast(1-x^\ast)$ reduces to $y^\ast>1-y^\ast$,
which holds because $y^\ast>\tfrac12$.

It remains to examine the other three combinations of regimes. If both
coordinates are active, then, as above, the equations
$\tfrac{\alpha x(1-y)}{z_p}=x$ and $\tfrac{\alpha y(1-x)}{z_r}=y$ give
$\alpha(1-y)=1$ and $\alpha(1-x)=1$, whence $x=y$. If both coordinates are
passive, the equations $\tfrac{1-x}{z_p}=x$ and $\tfrac{1-y}{z_r}=y$ give
$\bigl(1-x\bigr)^2=\alpha x^2(1-y)$ and $\bigl(1-y\bigr)^2=\alpha y^2(1-x)$; dividing the first
relation by the second, we get
$$\frac{\bigl(1-x\bigr)^3}{x^2}=\frac{\bigl(1-y\bigr)^3}{y^2},$$
and since the function $u\mapsto\bigl(1-u\bigr)^3/u^2$ is strictly decreasing on
$(0,1)$, again $x=y$. Finally, by the equivariance
$\Phi\circ\sigma=\sigma\circ\Phi$, a point $(x,y)$ is a fixed point with the
first coordinate passive and the second active if and only if $(y,x)$ is a
fixed point of the type considered above; hence the only such point is
$(y^\ast,x^\ast)$. Therefore the only off-diagonal fixed points are
$(x^\ast,y^\ast)$ and $(y^\ast,x^\ast)$.

In a neighbourhood of $(x^\ast,y^\ast)$ the first coordinate is given by the
active branch $x'=\alpha x(1-y)/z_p$ and the second by the passive branch
$y'=(1-y)/z_r$; both regime inequalities are strict, so $\Phi$ is smooth
there. By the quotient rule, the numerators of the partial derivatives
simplify identically: for $\partial x'/\partial x$ we have
$z_p-x\bigl(\alpha(1-y)-1\bigr)\equiv1$; for $\partial x'/\partial y$ the
numerator equals $-\alpha x\bigl(z_p-\alpha x(1-y)\bigr)=-\alpha x(1-x)$;
for $\partial y'/\partial y$ it equals
$-z_r-(1-y)\bigl(\alpha(1-x)-1\bigr)=-\alpha(1-x)$. Substituting
$z_p=1$, $z_r=\tfrac{1}{\alpha-1}$ and \eqref{eq_fixed}, we obtain
\begin{align*}
\frac{\partial x'}{\partial x}&=\frac{\alpha(1-y)}{z_p^2}=1, &
\frac{\partial x'}{\partial y}&=-\frac{\alpha x(1-x)}{z_p^2}
=-\frac{x^\ast}{\bigl(\alpha-1\bigr)^2},\\[1mm]
\frac{\partial y'}{\partial x}&=\frac{\alpha y(1-y)}{z_r^2}
=\frac{\bigl(\alpha-1\bigr)^3}{\alpha}, &
\frac{\partial y'}{\partial y}&=-\frac{\alpha(1-x)}{z_r^2}=-1 .
\end{align*}
Hence
$$J=\begin{pmatrix}
1 & -\dfrac{x^\ast}{\bigl(\alpha-1\bigr)^2}\\[3mm]
\dfrac{\bigl(\alpha-1\bigr)^3}{\alpha} & -1
\end{pmatrix},$$
$$ \operatorname{Tr}J=0,\qquad
\det J=-1+\frac{(\alpha-1)x^\ast}{\alpha}
=-\frac{\alpha^2-\alpha+1}{\alpha^2(\alpha-1)}=-\rho^2(\alpha).$$
Therefore the eigenvalues of $J$ are $\pm\rho(\alpha)$, and by the
Cayley--Hamilton theorem $J^2=(\operatorname{Tr}J)\,J-(\det J)\,I=\rho^2(\alpha)I$.

The inequality $\rho(\alpha)<1$ is equivalent to
$f(\alpha):=\alpha^3-2\alpha^2+\alpha-1>0$; since
$f'(\alpha)=(3\alpha-1)(\alpha-1)>0$ for $\alpha>1$ and $f(2)=1$, it holds
for all $\alpha>2$. Moreover,
$$\frac{d\rho^2}{d\alpha}=-\frac{\alpha^3-2\alpha^2+4\alpha-2}{\alpha^3\bigl(\alpha-1\bigr)^2}<0
\qquad(\alpha>2),$$
since $\alpha^3-2\alpha^2=\alpha^2(\alpha-2)>0$ and $4\alpha-2>0$; finally,
$\rho^2(2)=\tfrac34$ and $\rho^2(\alpha)\sim\tfrac1\alpha$ as $\alpha\to\infty$.
By the theorem on stability in the first approximation
\cite[Theorem~4.33]{elaydi}, the fixed point $(x^\ast,y^\ast)$ is locally
asymptotically stable, and since $J^2=\rho^2 I$, the deviation from it is
multiplied asymptotically by $\rho^2(\alpha)$ over every two steps. By
equivariance, the same holds for $(y^\ast,x^\ast)$.

Along a trajectory converging to $(x^\ast,y^\ast)$, the regime inequalities
established above hold from some step on, since they are strict at the
fixed point. From this step on, the maximum of $\mathbf p^t$ keeps its
position and its value tends to $x^\ast$, while the maximum of
$\mathbf r^t$ switches its position at every step and its value tends to
$y^\ast$.
\end{proof}

By Lemma~\ref{lem_list_2cycles}(2), for $\alpha>2$ the map $\Phi$ has no exchange $2$-cycles, so that the $4$-cycle of Theorem~\ref{thm_cycle_general} disappears and its role is taken over by the fixed points of Theorem~\ref{thm_fixed_general}. As for $0<\alpha<2$, numerical experiments indicate that generic trajectories with $x^0\neq y^0$ converge to these fixed points, whereas a rigorous proof of global convergence remains open. By Lemma~\ref{lem_order}, along every such trajectory the order of the maxima is preserved and $|\ell(x^t)-\ell(y^t)|\ge2\ln\tfrac\alpha2$ for all sufficiently large $t$, so that the trajectory does not approach the symmetric state.

\begin{remark}\label{rem_winner}
In applied terms, Theorem~\ref{thm_fixed_general} means that for $\alpha>2$ the model deterministically separates the players. The player whose maximum tends to the larger value $x^\ast$, the \emph{winner}, eventually keeps its dominant position and, as $\alpha$ grows, concentrates an increasing share of its mass there, since $x^\ast\to1$; the player with the value $y^\ast$, the \emph{loser}, remains in a perpetual switching regime. By Lemma~\ref{lem_order}(2), along a trajectory converging to one of these fixed points the winner is the player with the larger initial maximum, and its maximum keeps its position at every step. These terms are used in this sense below, in particular in the figure captions.
\end{remark}

\section{Bifurcation analysis}\label{bif}

The object of the analysis is the autonomous map of maxima $\Phi$ defined by
system~\eqref{eq_macro_system}, with $\alpha>0$ as the bifurcation parameter.
We show that the invariant sets found in Section~\ref{map} undergo a
bifurcation only at $\alpha=2$, and describe the nature of this transition.

\subsection{The symmetric branch and its saddle structure}

\begin{lem}\label{lem_diag}
Let $\alpha>0$. Then
\begin{enumerate}
\item the diagonal $\{x=y\}$ is invariant under $\Phi$, and the restriction of $\Phi$ to it is the map
\begin{equation}\label{eq_g}
g(x)=\frac{\max\{\alpha x,\,1\}}{\alpha x+1},
\end{equation}
which is piecewise smooth with a kink at $x=\tfrac1\alpha$;
\item for $\alpha\neq2$ the map $g$ has a unique fixed point $x_{\mathrm{sym}}(\alpha)$ in $\left(\tfrac12,1\right)$,
\begin{equation}\label{eq_xsym}
x_{\mathrm{sym}}(\alpha)=
\begin{cases}
\dfrac{-1+\sqrt{1+4\alpha}}{2\alpha}, & 0<\alpha<2,\\[2mm]
1-\dfrac1\alpha, & \alpha>2,
\end{cases}
\end{equation}
and $\alpha x_{\mathrm{sym}}<1$ for $0<\alpha<2$, $\alpha x_{\mathrm{sym}}>1$ for $\alpha>2$;
\item the function $x_{\mathrm{sym}}$, extended by $x_{\mathrm{sym}}(2)=\tfrac12$, is continuous on $(0,\infty)$, strictly decreasing on $(0,2]$ and strictly increasing on $[2,\infty)$;
\item $g'(x_{\mathrm{sym}})=x_{\mathrm{sym}}-1$ for $0<\alpha<2$ and $g'(x_{\mathrm{sym}})=\tfrac1\alpha$ for $\alpha>2$; in particular, $|g'(x_{\mathrm{sym}})|<1$ for $\alpha\neq2$, and for $\alpha=2$ the one-sided derivatives of $g$ at $x=\tfrac12$ are equal to $-\tfrac12$ and $\tfrac12$.
\end{enumerate}
\end{lem}

\begin{proof}
The invariance of the diagonal and the form of $g$ were established in the proof of Proposition~\ref{prop_sym_general}; the two arguments of the maximum in \eqref{eq_g} coincide at $\alpha x=1$, which gives the kink at $x=\tfrac1\alpha$.

We solve $g(x)=x$ on each branch. If $\alpha x<1$, then $g(x)=\tfrac{1}{\alpha x+1}$, and the equation $x(\alpha x+1)=1$ gives
$\alpha x^2+x-1=0$, whose positive root is $\tfrac{-1+\sqrt{1+4\alpha}}{2\alpha}$; the branch condition $\alpha x<1$
for this root reads $\tfrac{-1+\sqrt{1+4\alpha}}{2}<1$, i.e.\
$\sqrt{1+4\alpha}<3$, which is equivalent to $\alpha<2$. In this case the
root lies in $\left(\tfrac12,1\right)$, since the polynomial $\alpha x^2+x-1$ equals
$\tfrac\alpha4-\tfrac12<0$ at $x=\tfrac12$ and $\alpha>0$ at $x=1$. If
$\alpha x>1$, then $g(x)=\tfrac{\alpha x}{\alpha x+1}$, and the equation
$x(\alpha x+1)=\alpha x$ gives $x=1-\tfrac1\alpha$; the branch condition
$\alpha x=\alpha-1>1$ is equivalent to $\alpha>2$, and then
$1-\tfrac1\alpha\in\left(\tfrac12,1\right)$.

For $\alpha=2$ both formulas in \eqref{eq_xsym} give $\tfrac12$, so $x_{\mathrm{sym}}$ is continuous. On $(0,2)$ the value $x_{\mathrm{sym}}$ satisfies $F(x,\alpha):=\alpha x^2+x-1=0$, and by the implicit function theorem
$$\frac{dx_{\mathrm{sym}}}{d\alpha}=-\frac{F_\alpha}{F_x}=-\frac{x_{\mathrm{sym}}^2}{2\alpha x_{\mathrm{sym}}+1}<0,$$
so $x_{\mathrm{sym}}$ decreases strictly from $1$ as $\alpha\to0^+$ to $\tfrac12$ as $\alpha\to2^-$. On $(2,\infty)$ we have $\tfrac{d}{d\alpha}\left(1-\tfrac1\alpha\right)=\tfrac{1}{\alpha^2}>0$, so $x_{\mathrm{sym}}$ increases strictly from $\tfrac12$ to $1$. Hence the graph of $x_{\mathrm{sym}}$ is V-shaped with the minimum $x_{\mathrm{sym}}(2)=\tfrac12$.

On the switching branch $g'(x)=-\tfrac{\alpha}{\bigl(\alpha x+1\bigr)^2}$; at the point $x_{\mathrm{sym}}$ we have $\alpha x_{\mathrm{sym}}+1=\tfrac{1}{x_{\mathrm{sym}}}$ and
$\alpha x_{\mathrm{sym}}^2=1-x_{\mathrm{sym}}$, whence $g'(x_{\mathrm{sym}})=-\alpha x_{\mathrm{sym}}^2=x_{\mathrm{sym}}-1\in\left(-\tfrac12,0\right)$. On the preserving branch $g'(x)=\tfrac{\alpha}{\bigl(\alpha x+1\bigr)^2}$, and at the point $1-\tfrac1\alpha$ the denominator equals $\alpha^2$, so that $g'(x_{\mathrm{sym}})=\tfrac1\alpha<\tfrac12$. For $\alpha=2$ the two formulas give $\mp\tfrac{2}{\left(2\cdot\frac12+1\right)^2}=\mp\tfrac12$ at $x=\tfrac12$.
\end{proof}

\begin{prop}\label{prop_saddle}
Let $\alpha\neq2$ and $s=x_{\mathrm{sym}}(\alpha)$. Then
\begin{enumerate}
\item the map $\Phi$ is smooth in a neighbourhood of $(s,s)$, and the
vectors $\bigl(1,1\bigr)^\top$ and $\bigl(1,-1\bigr)^\top$ are eigenvectors of $J_\Phi(s,s)$
with eigenvalues
$$\lambda_\parallel=g'(s),\qquad
\lambda_\perp=\begin{cases}-(1+s),&0<\alpha<2,\\[1mm]\ \ \,1+s,&\alpha>2;\end{cases}$$
\item $|\lambda_\parallel|<1$ and $|\lambda_\perp|>\tfrac32$, so that
$(s,s)$ is a saddle of $\Phi$;
\item for $\alpha=1$ one has $\lambda_\perp=-\varphi$, where
$\varphi=\tfrac{1+\sqrt5}{2}$ is the golden ratio.
\end{enumerate}
\end{prop}

\begin{proof}
Let $0<\alpha<2$. At the point $(s,s)$ the strict inequality $\alpha s<1$ holds, so in a neighbourhood of this point both coordinates of $\Phi$ are given by the passive branches,
$$x'=\frac{1-x}{z_p},\quad z_p=\alpha x(1-y)+(1-x);\qquad y'=\frac{1-y}{z_r},\quad z_r=\alpha y(1-x)+(1-y),$$
and $\Phi$ is smooth there. By Lemma~\ref{lem_diag}, $\alpha s^2=1-s$, whence $\alpha s+1=\tfrac1s$ and $z_p(s,s)=z_r(s,s)=(1-s)(\alpha s+1)=\tfrac{1-s}{s}=:z$.

By the quotient rule,
$$\frac{\partial x'}{\partial x} =\frac{-z_p-(1-x)\bigl[\alpha(1-y)-1\bigr]}{z_p^2}.$$
Substituting $x=y=s$ and $\alpha(1-s)=\alpha\cdot\alpha s^2=\alpha^2s^2$, we obtain
$$\frac{\partial x'}{\partial x} = \frac{-\tfrac{1-s}{s}-(1-s)\bigl(\alpha^2s^2-1\bigr)}{\left(\tfrac{1-s}{s}\right)^{2}} = -\frac{s+\alpha^2s^4-s^2}{1-s}.$$
Since $\alpha^2s^4=\bigl(\alpha s^2\bigr)^2=\bigl(1-s\bigr)^2$, the numerator equals $s-s^2+\bigl(1-s\bigr)^2=1-s$, so that $\dfrac{\partial x'}{\partial x}=-1$. Further,
$$\frac{\partial x'}{\partial y}
=-\frac{(1-x)\,\partial z_p/\partial y}{z_p^2}
=\frac{(1-s)\,\alpha s}{z^{2}}
=\frac{\alpha s(1-s)\,s^2}{\bigl(1-s\bigr)^2}
=\frac{\alpha s^3}{1-s}=\frac{\alpha s^3}{\alpha s^2}=s.$$
By the symmetry $x\leftrightarrow y$, also $\dfrac{\partial y'}{\partial y}=-1$ and $\dfrac{\partial y'}{\partial x}=s$. Hence
$$J_\Phi(s,s)=\begin{pmatrix}-1&s\\ s&-1\end{pmatrix},$$
$$J_\Phi(s,s)\begin{pmatrix}1\\1\end{pmatrix}=(s-1)\begin{pmatrix}1\\1\end{pmatrix},\quad
J_\Phi(s,s)\begin{pmatrix}1\\-1\end{pmatrix}=-(1+s)\begin{pmatrix}1\\-1\end{pmatrix}.$$
By Lemma~\ref{lem_diag}(4), $\lambda_\parallel=s-1=g'(s)$ and $|\lambda_\parallel|<\tfrac12$, while $|\lambda_\perp|=1+s>\tfrac32$.

Let $\alpha>2$. Then $s=1-\tfrac1\alpha$ and $\alpha s=\alpha-1>1$, so in a
neighbourhood of $(s,s)$ both coordinates of $\Phi$ are given by the active
branches $x'=\alpha x(1-y)/z_p$ and $y'=\alpha y(1-x)/z_r$. We use
$1-s=\tfrac1\alpha$, $\alpha(1-s)=1$ and
$z_p(s,s)=(1-s)(\alpha s+1)=\tfrac1\alpha\cdot\alpha=1$. Since identically
$$z_p-x\bigl(\alpha(1-y)-1\bigr)=\alpha x(1-y)+1-x-\alpha x(1-y)+x=1,$$
we have
$$\frac{\partial x'}{\partial x}
=\frac{\alpha(1-y)\bigl(z_p-x(\alpha(1-y)-1)\bigr)}{z_p^2}
=\frac{\alpha(1-y)}{z_p^2}\ \Bigr|_{(s,s)}=\alpha(1-s)=1,$$
$$\frac{\partial x'}{\partial y}
=-\frac{\alpha x(1-x)}{z_p^2}\ \Bigr|_{(s,s)}=-\alpha s(1-s)=-s.$$
By symmetry,
$J_\Phi(s,s)=\bigl(\begin{smallmatrix}1&-s\\-s&1\end{smallmatrix}\bigr)$,
with the eigenvalue $\lambda_\parallel=1-s=\tfrac1\alpha=g'(s)$ on
$\bigl(1,1\bigr)^\top$ and $\lambda_\perp=1+s=2-\tfrac1\alpha>\tfrac32$ on
$\bigl(1,-1\bigr)^\top$.

For $\alpha=1$ we have $s=\varphi^{-1}=\tfrac{\sqrt5-1}{2}$ and
$1+s=1+\varphi^{-1}=\varphi$ by the defining property $\varphi=1+\tfrac1\varphi$
of the golden ratio, so that $\lambda_\perp=-\varphi\approx-1.618$.
\end{proof}

\begin{remark}\label{rem_saddle_set}
By the stable manifold theorem \cite[Ch.~2]{palis1982}, the points of a neighbourhood of $(s,s)$
whose trajectories converge to $(s,s)$ form a smooth curve tangent to
$\bigl(1,1\bigr)^\top$. Since the diagonal is invariant and, by
Lemma~\ref{lem_diag}(4), attracting along itself, this curve is a segment of
the diagonal. Consequently, every trajectory of $\Phi$ converging to $(s,s)$
lies on the diagonal from some step on, i.e.\ the set of initial states
attracted by the symmetric state is contained in
$\bigcup_{t\ge0}\Phi^{-t}\bigl(\{x=y\}\bigr)$; for $\alpha\in(1,2)$ this
set contains the curve $S$ of Remark~\ref{rem_exceptional}. On the branches
of $\Phi$ of the same type, $x'-y'=\pm\alpha(x-y)(1-xy)/(z_pz_r)$, so that a
point off the diagonal and off $T_p\cup T_r$ is mapped onto the diagonal only
if it lies on $S$. Hence
$\bigcup_{t\ge0}\Phi^{-t}\bigl(\{x=y\}\bigr)\subset\{x=y\}\cup E_\alpha$, and
this set has Lebesgue measure zero. In a neighbourhood of $(s,s)$ the
distance of a trajectory to the diagonal grows by the factor
$|\lambda_\perp|=1+x_{\mathrm{sym}}(\alpha)$ at every step.
\end{remark}

\begin{remark}\label{rem_multipliers}
The multipliers in Theorems~\ref{thm_cycle_general} and~\ref{thm_fixed_general}
and in Proposition~\ref{prop_saddle} can also be read off from
Lemma~\ref{lem_logform}. The map $(x,y)\mapsto(\ell(x),\ell(y))$ is a
diffeomorphism of $\left(\tfrac12,1\right)^2$ onto $\bigl(0,\infty\bigr)^2$ that conjugates $\Phi$
and $\Psi$, so that corresponding periodic points have the same multipliers,
and the Jacobian matrix of $\Psi$ depends on $\alpha$ only through the types
of the branches. As shown in the proofs of Theorems~\ref{thm_cycle_general}
and~\ref{thm_fixed_general}, at the point $(x^\ast,y^\ast)$ of the exchange
$2$-cycle and at the fixed point $(x^\ast,y^\ast)$ the first coordinate of
$\Phi$ is given by the active branch and the second by the passive one. For
the exchange $2$-cycle, the monodromy matrix of $\Psi$ is therefore
$$\begin{pmatrix}-1&x^\ast\\-y^\ast&1\end{pmatrix}
\begin{pmatrix}1&-y^\ast\\x^\ast&-1\end{pmatrix};$$
its determinant equals $\bigl(1-x^\ast y^\ast\bigr)^2$, and its trace
$\bigl(x^\ast\bigr)^2+\bigl(y^\ast\bigr)^2-2$ has modulus less than $2(1-x^\ast y^\ast)$, so that
the multipliers are complex conjugate with modulus $1-x^\ast y^\ast$. For the
fixed point, the matrix
$\bigl(\begin{smallmatrix}1&-y^\ast\\x^\ast&-1\end{smallmatrix}\bigr)$ has
zero trace and determinant $-(1-x^\ast y^\ast)$, so that its eigenvalues are
$\pm\sqrt{1-x^\ast y^\ast}=\pm\rho(\alpha)$. At the symmetric state $(s,s)$
both coordinates are given by the passive branch for $\alpha<2$ and by the
active branch for $\alpha>2$; accordingly, the Jacobian matrix of $\Psi$
equals $\mp\bigl(\begin{smallmatrix}1&-s\\-s&1\end{smallmatrix}\bigr)$, with
the eigenvalues $\mp(1-s)$ along the diagonal and $\mp(1+s)$ transversally to
it, and the change of sign of $\lambda_\perp$ at $\alpha=2$ is the change of
the types of both branches. Finally, the curves $\{\eta(u,v)=0\}$ and $\{\eta(v,u)=0\}$ have Lebesgue
measure zero, and so do their images under the locally Lipschitz map $\Psi$;
on each of the open sets between them the types of the branches are fixed,
and by Lemma~\ref{lem_logform}(4) $\Psi$ is injective there with
$|\det J_\Psi|\le\tfrac34$. Hence
$\operatorname{Leb}(\Psi(B))\le\tfrac34\operatorname{Leb}(B)$ for every
measurable set $B\subset[0,\infty)^2$. Let $I$ be a measurable set with
$\Phi(I)=I$, for instance an $\omega$-limit set of $\Phi$. Then
$I=\Phi^2(I)\subset K_\alpha$, since every trajectory enters $K_\alpha$ after
two steps, so that the image $B$ of $I$ under $(x,y)\mapsto(\ell(x),\ell(y))$
is bounded and satisfies $\Psi(B)=B$; hence $\operatorname{Leb}(B)=0$, and
$I$ has Lebesgue measure zero.
\end{remark}

\subsection{The transition at $\alpha=2$}

\begin{theorem}\label{thm_bifurcation}
Let $x_{\mathrm{sym}}(\alpha)$ be given by \eqref{eq_xsym}. Then
\begin{enumerate}
\item for $0<\alpha<2$ the fixed point $x_{\mathrm{sym}}$ of $g$ lies on the
switching branch, $x_{\mathrm{sym}}<\tfrac1\alpha$, and for $\alpha>2$ on the
preserving branch, $x_{\mathrm{sym}}>\tfrac1\alpha$; the equality
$x_{\mathrm{sym}}(\alpha)=\tfrac1\alpha$ holds only for $\alpha=2$;
\item the map $\Phi$ has an exchange $2$-cycle if and only if $0<\alpha<2$,
and fixed points off the diagonal if and only if $\alpha>2$; the points
$(x^\ast,y^\ast)$ of both the exchange $2$-cycle of
Theorem~\ref{thm_cycle_general} and the fixed point of
Theorem~\ref{thm_fixed_general} tend to $\left(\tfrac12,\tfrac12\right)$ as
$\alpha\to2^-$ and $\alpha\to2^+$, respectively;
\item for $\alpha\in(0,2)\cup(2,\infty)$ the symmetric fixed point
$(x_{\mathrm{sym}},x_{\mathrm{sym}})$, the exchange $2$-cycle and the
off-diagonal fixed points depend analytically on $\alpha$, satisfy their
regime inequalities strictly, and none of their multipliers lies on the
unit circle;
\item the transverse multiplier $\lambda_\perp$ of the symmetric fixed point
tends to $-\tfrac32$ as $\alpha\to2^-$ and to $\tfrac32$ as $\alpha\to2^+$.
\end{enumerate}
\end{theorem}

\begin{proof}
By Lemma~\ref{lem_diag}(2), $\alpha x_{\mathrm{sym}}<1$ for $0<\alpha<2$ and $\alpha x_{\mathrm{sym}}>1$ for $\alpha>2$, while $x_{\mathrm{sym}}(2)=\tfrac12=\tfrac1\alpha$ for $\alpha=2$. On the switching branch the position of the maximum changes at every step, which gives the synchronous and mirror $2$-cycles of Proposition~\ref{prop_sym_general}; on the preserving branch the positions of the maxima are constant, which gives the fixed points of Proposition~\ref{prop_sym_fixed}.

By Lemma~\ref{lem_list_2cycles}, exchange $2$-cycles exist only for $0<\alpha<2$, and by Remark~\ref{rem_alpha2_cycle} their values tend to $\tfrac12$ as $\alpha\to2^-$. In the proof of Theorem~\ref{thm_fixed_general} the condition $\alpha>2$ is used only to check admissibility; hence for every $\alpha>0$ the only possible off-diagonal fixed points of $\Phi$ are $(x^\ast,y^\ast)$ and $(y^\ast,x^\ast)$ with $x^\ast,y^\ast$ given by \eqref{eq_fixed}. For $\alpha\le2$ we have $y^\ast=1-\tfrac1\alpha\le\tfrac12$, so these points do not belong to $\left(\tfrac12,1\right)^2$. As $\alpha\to2^+$,
$$y^\ast=1-\frac1\alpha\to\frac12,\qquad x^\ast=1-\frac{1}{\alpha\bigl(\alpha-1\bigr)^2}\to1-\frac{1}{2\cdot1}=\frac12.$$

The function $x_{\mathrm{sym}}$ is given explicitly by \eqref{eq_xsym}; the value $x^\ast$ of the exchange $2$-cycle is a simple root of the cubic \eqref{eq_cubic_general}, since $P'(x^\ast)>0$, and therefore depends analytically on $\alpha$ by the implicit function theorem; the fixed point \eqref{eq_fixed} is given explicitly. The regime inequalities are strict by Lemma~\ref{lem_diag}(2) and by the proofs of Theorems~\ref{thm_cycle_general} and~\ref{thm_fixed_general}. The multipliers of the exchange $2$-cycle have modulus $1-x^\ast y^\ast<\tfrac34$ (Theorem~\ref{thm_cycle_general}), the eigenvalues $\pm\rho(\alpha)$ at the off-diagonal fixed points satisfy $\rho(\alpha)<\tfrac{\sqrt3}{2}$ (Theorem~\ref{thm_fixed_general}), and the multipliers of the symmetric fixed point satisfy $|\lambda_\parallel|<1<\tfrac32<|\lambda_\perp|$ (Proposition~\ref{prop_saddle}). Finally, by Proposition~\ref{prop_saddle}, the transverse multiplier of the symmetric fixed point equals $\lambda_\perp=-(1+x_{\mathrm{sym}})$ for $0<\alpha<2$ and $\lambda_\perp=1+x_{\mathrm{sym}}$ for $\alpha>2$; since $x_{\mathrm{sym}}$ is continuous and $x_{\mathrm{sym}}(2)=\tfrac12$ (Lemma~\ref{lem_diag}), it tends to $-\tfrac32$ as $\alpha\to2^-$ and to $\tfrac32$ as $\alpha\to2^+$.
\end{proof}

\begin{figure}[htbp]
\centering
\includegraphics[width=\textwidth]{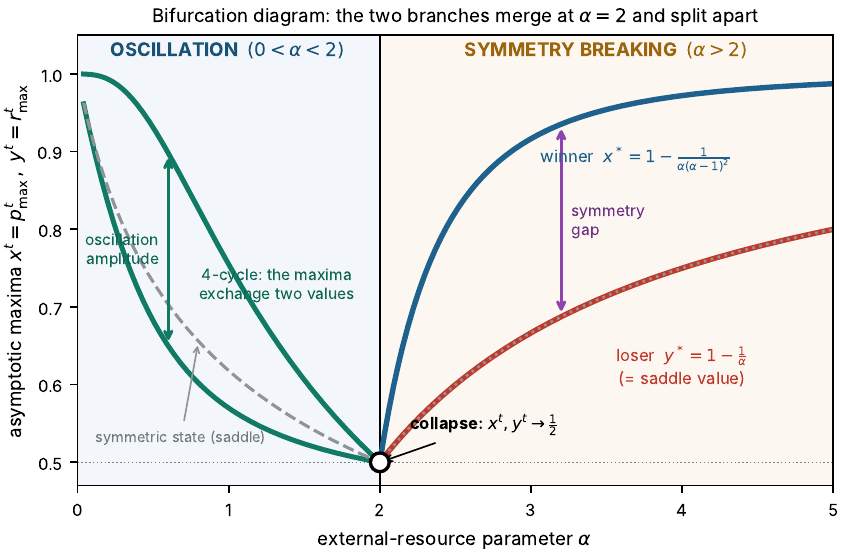}
\caption{Asymptotic values of the maxima $x^t,y^t$ as functions of $\alpha$.
For $0<\alpha<2$ (teal): the two values $x^\ast>y^\ast$ of the $4$-cycle
(Theorem~\ref{thm_cycle_general}), which both maxima take in turn; for
$\alpha>2$: the coordinates $x^\ast$ (blue) and $y^\ast$ (red) of the
asymmetric fixed point (Theorem~\ref{thm_fixed_general}). The grey curve
is the symmetric state $x_{\mathrm{sym}}(\alpha)$, which is a saddle of $\Phi$
for $\alpha\neq2$; for $\alpha>2$ it coincides with $y^\ast=1-\tfrac1\alpha$
(dotted).
The open circle marks the point $\left(2,\tfrac12\right)$: at $\alpha=2$ both maxima
tend to $\tfrac12$ (Theorem~\ref{thm_alpha2}).}
\label{fig_bif_main}
\end{figure}

Theorem~\ref{thm_bifurcation} shows that $\alpha=2$ is the only value of the parameter at which the invariant sets of $\Phi$ found in Section~\ref{map} bifurcate. At this value every trajectory converges to $\left(\tfrac12,\tfrac12\right)$ (Theorem~\ref{thm_alpha2}). On the invariant diagonal the transition is a border collision in the sense of \cite{nusse1992,dibernardo2008,simpson2016}: the fixed point of the piecewise-smooth map $g$ crosses its kink $x=\tfrac1\alpha$, and in the full system the period of the symmetric regime changes from $2$ to $1$. The kink corresponds to the point $\left(\tfrac1\alpha,\tfrac1\alpha\right)$, at which the diagonal meets both switching curves $T_p$ and $T_r$, and for $\alpha=2$ this point is the corner $\left(\tfrac12,\tfrac12\right)$ of the phase space. Thus, unlike a generic border collision, in which a fixed point crosses a single switching manifold, here the symmetric fixed point reaches the intersection of the two switching curves, and it does so at the boundary of the phase space. In the plane this is accompanied by symmetry breaking: the off-diagonal fixed points branch off from the corner $\left(\tfrac12,\tfrac12\right)$, into which the exchange $2$-cycle shrinks from the other side. For $\alpha>2$ the distance between the coordinates of the asymmetric fixed point, $x^\ast-y^\ast=\frac1\alpha\left(1-\frac{1}{\bigl(\alpha-1\bigr)^2}\right)$, is positive and tends to $0$ as $\alpha\to2^+$; together with Remark~\ref{rem_alpha2_cycle} and Theorem~\ref{thm_fixed_general} this describes quantitatively the death of the oscillations and the opening of the symmetry gap (Fig.~\ref{fig_bif_main}).

\begin{figure}[htbp]
\centering
\includegraphics[width=\textwidth]{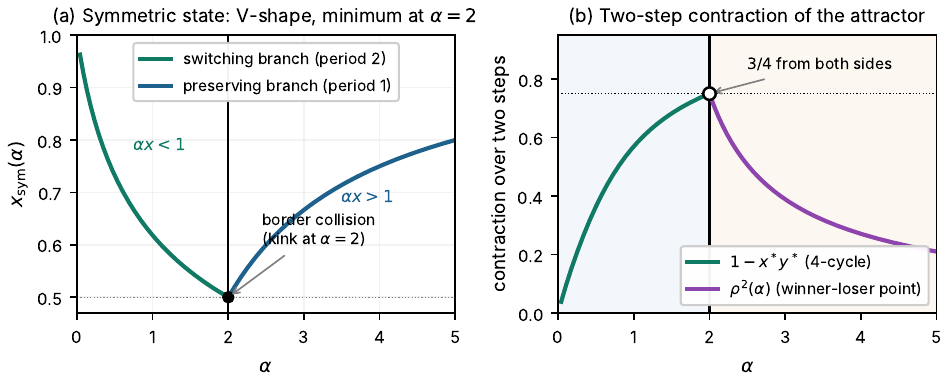}
\caption{(a) The symmetric state $x_{\mathrm{sym}}(\alpha)$: the fixed point
of the diagonal map $g$ lies on the switching branch $\alpha x<1$ for
$\alpha<2$ (period~$2$ in the full system) and on the preserving branch
$\alpha x>1$ for $\alpha>2$ (period~$1$); it crosses the kink of $g$ at
$\alpha=2$. (b) Contraction of the asymmetric attractor over two steps:
the modulus $1-x^\ast y^\ast$ of the multipliers of the $4$-cycle for
$\alpha<2$ and $\rho^2(\alpha)$ for $\alpha>2$ (recall that $J^2=\rho^2 I$).
Both tend to $\tfrac34$ as $\alpha\to2$.}
\label{fig_bif_mech}
\end{figure}

\begin{figure}[htbp]
\centering
\includegraphics[width=\textwidth]{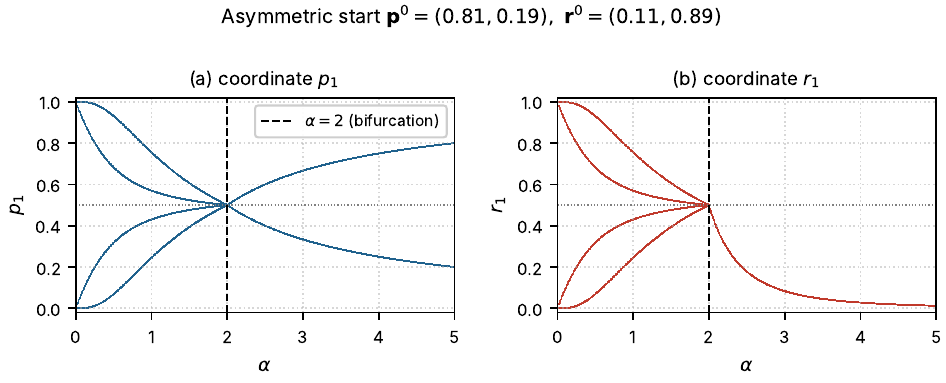}
\caption{Coordinate-wise bifurcation diagram of the full system for the
asymmetric initial state $\mathbf p^0=(0.81,0.19)$, $\mathbf r^0=(0.11,0.89)$:
the last $8$ of $3008$ iterations for $600$ values of $\alpha$; the second
coordinates are $p_2=1-p_1$, $r_2=1-r_1$. (a) $p_1$; (b) $r_1$. For
$0<\alpha<2$ each coordinate takes four values (the $4$-cycle). For
$\alpha>2$ the vector $\mathbf r$, which has the larger initial maximum,
converges to a stationary distribution with $r_1\to1-x^\ast$, while $p_1$
alternates between $y^\ast$ and $1-y^\ast$.}
\label{fig_bif_coord}
\end{figure}

\begin{figure}[htbp]
\centering
\includegraphics[width=\textwidth]{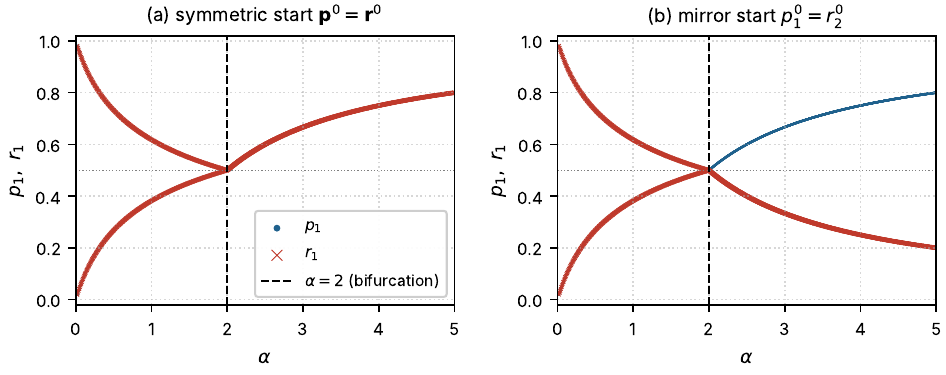}
\caption{The same diagram for symmetric initial states:
(a) $\mathbf p^0=\mathbf r^0=(0.81,0.19)$, where $p_1^t=r_1^t$ for all $t$;
(b) $\mathbf p^0=(0.81,0.19)$, $\mathbf r^0=(0.19,0.81)$, where
$r_1^t=1-p_1^t$. For $0<\alpha<2$ the system converges to the synchronous,
respectively mirror, $2$-cycle with value $x_{\mathrm{sym}}(\alpha)$
(Proposition~\ref{prop_sym_general}); for $\alpha>2$ to the fixed points of
Proposition~\ref{prop_sym_fixed}.}
\label{fig_bif_coord_sym_mirror}
\end{figure}

Fig.~\ref{fig_bif_main} summarizes the qualitative picture. As $\alpha$
increases to $2$, the two branches approach each other, merge at the point
$\left(2,\tfrac12\right)$ and separate again for $\alpha>2$. The meaning of the two
branches changes at $\alpha=2$. For $0<\alpha<2$ they are the two values
$x^\ast$ and $y^\ast$ which the maximum of each vector takes alternately
along the $4$-cycle (Theorem~\ref{thm_cycle_general}); for this reason both
branches are drawn in the same colour. For $\alpha>2$ they are the limits of
the maxima of two different vectors: the maximum of one of them tends to
$x^\ast$, and $x^\ast\to1$ as $\alpha\to\infty$, while the maximum of the
other one tends to $y^\ast=1-\tfrac1\alpha$
(Theorem~\ref{thm_fixed_general}). Thus at $\alpha=2$ the alternation of
values in time is replaced by a permanent asymmetry between the vectors.

Fig.~\ref{fig_bif_mech} illustrates the mechanism of the transition.
Panel~(a) shows the V-shaped graph of $x_{\mathrm{sym}}$
(Lemma~\ref{lem_diag}): the symmetric fixed point exists for all
$\alpha\neq2$, but at $\alpha=2$, where $x_{\mathrm{sym}}=\tfrac12$, it
passes through the kink $x=\tfrac1\alpha$ of $g$ and moves from the
switching branch to the preserving one. Panel~(b) shows the contraction of
the asymmetric attractors over two steps: the modulus $1-x^\ast y^\ast$ of
the multipliers of the $4$-cycle for $\alpha<2$ and $\rho^2(\alpha)$ for
$\alpha>2$. Both quantities tend to $\tfrac34$ as $\alpha\to2$
(Remark~\ref{rem_alpha2_cycle} and Theorem~\ref{thm_fixed_general}), so that the convergence to the
asymmetric attractor slows down as $\alpha$ approaches the bifurcation value
from either side.

Fig.~\ref{fig_bif_coord} shows the coordinates $p_1$ and $r_1$ of the full
system and thereby reveals the structure hidden in the diagram of the
maxima. For $0<\alpha<2$ each coordinate takes the four values
$x^\ast$, $y^\ast$, $1-x^\ast$, $1-y^\ast$ along the cycle, which makes the
period $4$ visible. For $\alpha>2$ the panel of $r_1$ shows a single branch,
since the vector $\mathbf r$ converges to a stationary distribution, whereas
the panel of $p_1$ shows two branches, since the maximum of $\mathbf p$
changes its position at every step (Theorem~\ref{thm_fixed_general}(3)).
Here $r^0_{\max}=0.89>p^0_{\max}=0.81$, and
the winner is $\mathbf r$, in agreement with Remark~\ref{rem_winner}.

To clarify the role of the initial symmetry, Fig.~\ref{fig_bif_coord_sym_mirror}
shows the same diagram for symmetric and mirror initial states. By
Propositions~\ref{prop_sym_general} and~\ref{prop_sym_fixed}, for
$\mathbf p^0=\mathbf r^0$ (panel~(a)) the equality $\mathbf p^t=\mathbf r^t$
holds for all $t$ and all $\alpha$, and for $p^0_1=r^0_2$ (panel~(b)) the
equality $r^t_1=1-p^t_1$ holds. Compared with Fig.~\ref{fig_bif_coord}, this
leads to two differences:
\begin{itemize}
\item for $0<\alpha<2$ each coordinate takes only the two values
$x_{\mathrm{sym}}$ and $1-x_{\mathrm{sym}}$, which corresponds to the
synchronous, respectively mirror, $2$-cycle, instead of the four values of
the $4$-cycle;
\item for $\alpha>2$ no symmetry breaking occurs: both vectors converge to
stationary distributions with the common maximal value $1-\tfrac1\alpha$
(period~$1$); this state is attracting within the invariant diagonal but is
a saddle of $\Phi$ (Proposition~\ref{prop_saddle}).
\end{itemize}
Thus the asymmetric attractors can be reached only from initial states off
the invariant symmetric set, and by Remark~\ref{rem_saddle_set} the set of
initial states attracted by the symmetric state has measure zero.

\section{Numerical experiments}\label{numerics}

In this section we test numerically the global convergence to the
asymmetric attractors, for which Theorems~\ref{thm_cycle_general}
and~\ref{thm_fixed_general} establish only local stability, and quantify the collapse at $\alpha=2$.
All computations are performed in double precision for the map $\Phi$ of
\eqref{eq_macro_system}; the coordinate-wise diagrams of
Figs.~\ref{fig_bif_coord} and~\ref{fig_bif_coord_sym_mirror} are computed for the
full system \eqref{conf_max}, with the normalizing denominators evaluated as the
sums of the numerators, so that the vectors remain stochastic up to rounding errors.
The code, with fixed random seeds, is available at
\url{https://github.com/kseniajasko/dynamic-conflict-systems} (folder \texttt{maxmap}) in the form of a Python package: the command \texttt{maxmap figures} produces the figures, and
the command \texttt{maxmap checks -{}-extended} performs all the checks described in
this section.

\subsection{Global convergence and phase portraits}

For each of the twelve values
$\alpha\in\{0.2,\,0.5,\,1,\,1.5,\,1.9,\,1.99,\,2.01,\,2.1,\,2.5,\,3,\,5,\,10\}$
we took $10^5$ initial points uniformly distributed in $\left(\tfrac12,1\right)^2$ and
iterated $\Phi$ up to $6000$ steps. A trajectory was counted as convergent
at the first even step at which its distance to the asymmetric attractor
(the exchange $2$-cycle for $\alpha<2$, the pair of fixed points of
Theorem~\ref{thm_fixed_general} for $\alpha>2$) was less than $10^{-10}$.
All trajectories converged, each in at most $156$ steps. The median numbers
of steps are $28$, $44$, $76$, $106$, $126$, $126$ for $\alpha=0.2$, $0.5$,
$1$, $1.5$, $1.9$, $1.99$ and $146$, $112$, $64$, $46$, $30$, $22$ for
$\alpha=2.01$, $2.1$, $2.5$, $3$, $5$, $10$; they grow as $\alpha\to2$, in
agreement with the contraction factor over two steps, which tends to
$\tfrac34$ from both sides (Fig.~\ref{fig_bif_mech}(b)). For $\alpha>2$, in
every case the limit fixed point was $(x^\ast,y^\ast)$ if $x^0>y^0$ and
$(y^\ast,x^\ast)$ if $x^0<y^0$: the player with the larger initial maximum
becomes the winner in the sense of Remark~\ref{rem_winner}, as required by
Lemma~\ref{lem_order}(2). The value $\alpha=2$ was not included, since
there the asymmetric attractors do not exist and every trajectory converges
to $\left(\tfrac12,\tfrac12\right)$ (Theorem~\ref{thm_alpha2}).
With the random seeds $1,\dots,8$ all trajectories also converged, the medians
were the same, and the maximal numbers of steps differed from those above by at
most $6$. A Newton search for periodic orbits of $\Phi$ of periods $1,\dots,8$, from
$2\cdot10^4$ random starting points for each period and each
$\alpha\in\{0.5,1,1.5,1.9,1.99,2.5,3,5\}$, found no periodic orbits off the switching
curves other than the symmetric fixed point and, for $\alpha<2$, the exchange
$2$-cycle or, for $\alpha>2$, the fixed points of Theorem~\ref{thm_fixed_general}.

\begin{figure}[htbp]
\centering
\includegraphics[width=\textwidth]{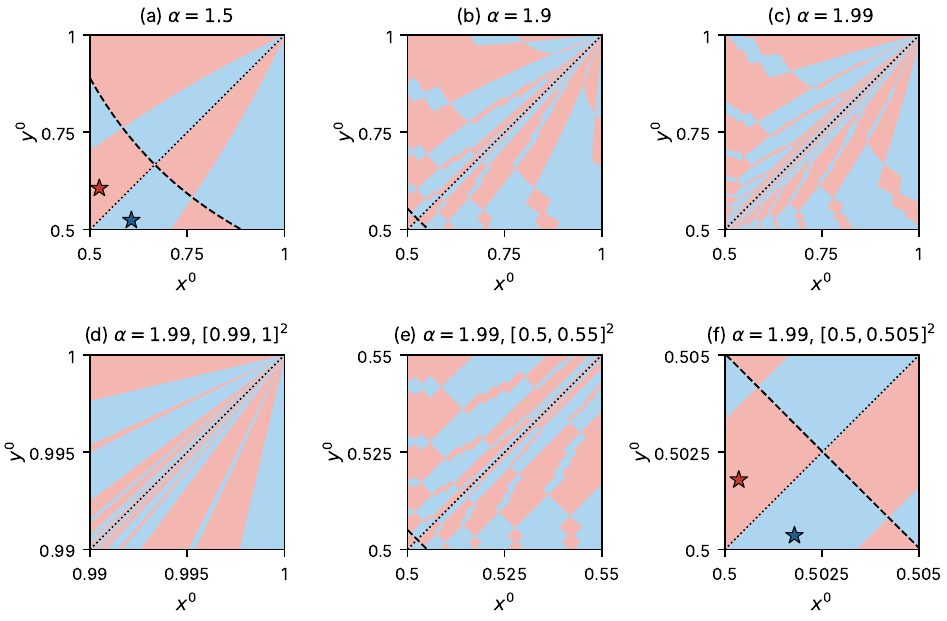}
\caption{Phases of the exchange $2$-cycle on grids of $600\times599$ initial
points $(x^0,y^0)$ after $4000$ iterations of $\Phi$: the colour shows which
point of the cycle the trajectory occupies at even times (light blue:
$(x^\ast,y^\ast)$, light red: $(y^\ast,x^\ast)$); all grid points converge to
the cycle. Stars: the points of the cycle (where they are resolved); dotted:
the diagonal; dashed: the curve $S$. (a)--(c) $\alpha=1.5$, $1.9$ and $1.99$
on the whole square. (d) $\alpha=1.99$, the window $[0.99,1]^2$ at the corner
$(1,1)$: the colour is constant along the rays through the corner. (e), (f)
$\alpha=1.99$, the windows $[0.5,0.55]^2$ and $[0.5,0.505]^2$ at the corner
$\left(\tfrac12,\tfrac12\right)$: the boundaries consist of segments of slope $\pm1$; in
(f) the diagonal and the curve $S$ cross at $\left(\tfrac1\alpha,\tfrac1\alpha\right)$.}
\label{fig_basins}
\end{figure}

The same conclusion holds on regular grids. Fig.~\ref{fig_basins} shows,
for grids of $600\times599$ initial points that contain no points of the
diagonal, which of the two points $(x^\ast,y^\ast)$, $(y^\ast,x^\ast)$ of the
exchange $2$-cycle the trajectory occupies at even times after $4000$
iterations; in all panels every grid point converges to the cycle. The two
colours are therefore not basins of different attractors but the two phases
of one attractor. Since $\Phi$ is continuous and the cycle is locally
attracting, the phase is locally constant on the set of initial points whose
trajectories converge to the cycle, and it can change only across points
whose trajectories do not converge to it, such as the points of the stable set
$$W=\bigcup_{t\ge0}\Phi^{-t}\bigl(\{x=y\}\bigr)$$
of the symmetric fixed point $(x_{\mathrm{sym}},x_{\mathrm{sym}})$, which is a
saddle by Proposition~\ref{prop_saddle}. For $200$ randomly chosen pairs of horizontally neighbouring grid points of
different colours in each panel, we located a point of the boundary by $50$
bisection steps, that is, to the resolution of double precision. Within $400$
steps, the trajectories of all $1200$ such points reach the diagonal to within
$10^{-13}$ and then approach the saddle to within $10^{-9}$, whereas the
trajectories of $10^4$ random initial points stay farther than $6\cdot10^{-5}$,
$4\cdot10^{-6}$ and $8\cdot10^{-7}$ from the saddle for $\alpha=1.5$, $1.9$ and
$1.99$, respectively; thus, up to rounding errors, the sampled boundary points
belong to $W$. A point off the diagonal is mapped onto the
diagonal only if it lies on the curve $S$ of Remark~\ref{rem_exceptional},
since for two branches of the same type the equality $x'=y'$ reduces to
$x\bigl(1-y\bigr)^2=y\bigl(1-x\bigr)^2$, i.e.\ to $x=y$. Hence $W$ consists of the diagonal and
the preimages of $S$. These preimages branch, since $\Phi$ is not invertible. By Lemma~\ref{lem_logform}, near a point of $T_p$ off $T_r$ the map $\Phi$ is conjugate to a diffeomorphism followed by the fold $(\xi,\zeta)\mapsto(|\xi|,\pm\zeta)$, so that such a point has a neighbourhood $U$ such that every point of $\Phi(U)$ off the side $x=\tfrac12$ has exactly two preimages in $U$, one on each side of $T_p$; the same holds for $T_r$. On the other hand, each point has at most one preimage in each of the regions A/A, A/P, P/A, P/P. For $0<\alpha\le1$ the curve $S$ is empty, and the
phase is determined by the sign of $x^0-y^0$; for $1<\alpha<2$ the preimages
of $S$ form a pattern that becomes increasingly intricate as $\alpha\to2$
(panels (a)--(c)). For $\alpha>2$ the two colours are the basins of the two
winner--loser fixed points, separated exactly by the diagonal by
Lemma~\ref{lem_order}(2).

The pattern is self-similar near both corners of the square. Near the corner
$(1,1)$ put $x=1-p$, $y=1-q$. The quantities $A$ and $B$ from the proof of
Lemma~\ref{lem_order} are $A=\alpha xq/p$ and $B=\alpha yp/q$, so that, up to
terms of order $p+q$, the image $\Phi(x,y)=\bigl(\chi(A),\chi(B)\bigr)$
depends only on the ratio $q/p$. Therefore, near $(1,1)$ the colour is
constant along the rays $q/p=\mathrm{const}$, and the pattern is invariant
under scaling about this corner: the windows $[0.9,1]^2$ and $[0.999,1]^2$, scaled onto panel~(d), differ
from it in about $5\%$ and $2\%$ of the points of the $600\times599$ grids;
about half of the differing points lie on rays that deviate from the diagonal by
less than $10^\circ$ (a sector containing $30\%$ of the grid points), where the
rays accumulate. Near the corner
$\left(\tfrac12,\tfrac12\right)$ put $x=\tfrac12+u$, $y=\tfrac12+v$. Up to terms of
second order in $u$, $v$ and $2-\alpha$, the map $\Phi$ reads
$$u'=\Bigl|u-\frac v2-\frac{2-\alpha}8\Bigr|,\qquad
v'=\Bigl|v-\frac u2-\frac{2-\alpha}8\Bigr|,$$
and in the variables $U=8u/(2-\alpha)$, $V=8v/(2-\alpha)$ it becomes the map
$$(U,V)\mapsto\left(\Bigl|U-\frac V2-1\Bigr|,\ \Bigl|V-\frac U2-1\Bigr|\right),$$
which contains no parameter. Consequently, at the scale $2-\alpha$ the pattern near the corner
is the same for all $\alpha$ close to $2$: for $\alpha=1.999$ the window
$[0.5,0.5005]^2$ differs from panel~(f) in less than $1\%$ of the grid
points. The exchange $2$-cycle of this map is
$\bigl\{\left(\tfrac{10}7,\tfrac27\right),\left(\tfrac27,\tfrac{10}7\right)\bigr\}$, in agreement
with Remark~\ref{rem_alpha2_cycle}. At the scales $2-\alpha\ll u+v\ll1$ the
constant terms are negligible, and the map is close to
$(u,v)\mapsto\left(|u-\tfrac v2|,\,|v-\tfrac u2|\right)$. On the cone
$0\le v\le\tfrac u2$ the latter acts as $(u,v)\mapsto\left(u-\tfrac v2,\tfrac u2-v\right)$,
maps this cone into itself, and its square is
$(u,v)\mapsto\left(\tfrac34u,\tfrac34v\right)$; the same holds on the cone
$0\le u\le\tfrac v2$. Since $w$ and $\Phi^2(w)$ always have the same colour, in this range of
scales the pattern approximately repeats under scaling by $\tfrac43$ about the
corner, and the number of such layers is of order $\log_{4/3}\frac1{2-\alpha}$.
Numerically, along the rays $v=su$, $s\in\{0.1,\,0.25,\,0.4\}$,
$10^{-6}\le u\le0.2$, the colour changes $17$--$18$, $32$--$33$ and $48$ times
for $\alpha=1.99$, $1.999$ and $1.9999$, that is, about $16=2\log_{4/3}10$ more
times for every factor $10$ in $\frac1{2-\alpha}$, twice per layer. The
repetition itself requires $2-\alpha\ll u+v\ll1$; for $\alpha=1.99$ this range
is too short, and the window of panel~(e) and its image under the scaling by
$\tfrac34$ about the corner differ in about $40\%$ of the grid points, whereas
for $\alpha=1.99999$ the window $[0.5,0.505]^2$ and its image differ in
$2.7\%$ of them. Finally, the preimages of lines of slope $\pm1$
under the affine pieces of the map above are again lines of slope $\pm1$;
since the diagonal and, near the corner, the curve $S$ have slopes $1$ and
$-1$, the colour boundaries near the corner consist of segments of slope
$\pm1$ (panels (e), (f)).

Although the pattern becomes increasingly intricate as $\alpha\to2$, it is
not fractal. The set $W$ lies in the union of the diagonal and the set
$E_\alpha$ of \eqref{eq_exceptional}, an at most countable union of smooth
curves, so that its Hausdorff dimension equals $1$; the self-similarity near
$\left(\tfrac12,\tfrac12\right)$ extends only over the finite range of scales between
$2-\alpha$ and $1$, and the one near $(1,1)$ is that of a family of rays.
Numerically, the box-counting dimension of the boundaries also equals $1$.
Let $f(\varepsilon)$ be the fraction of $10^6$ random initial points whose
colour changes under a displacement of length $\varepsilon$ in a random
direction; for boundaries of box-counting dimension $D$ one has
$f(\varepsilon)\sim\varepsilon^{2-D}$ \cite{mcdonald1985}. For $\varepsilon=10^{-4}$ and $10^{-5}$ we found
$f(\varepsilon)/\varepsilon\approx6.6$ and $7.5$, $24.8$ and $26.4$, $58.5$ and
$62.5$ for $\alpha=1.5$, $1.9$ and $1.99$, respectively; the corresponding
estimates of $2-D$ are $0.95\pm0.05$, $0.97\pm0.03$ and $0.97\pm0.02$ (binomial
standard errors), in agreement with $D=1$. For $\alpha$ close to $2$ this ratio is smaller at coarser scales
(for instance, $22$ at $\varepsilon=10^{-2}$ for $\alpha=1.99$), which
reflects the layered structure near the corner $\left(\tfrac12,\tfrac12\right)$.

\begin{remark}\label{rem_global_open}
These experiments indicate that for $0<\alpha<2$ almost every trajectory
converges to the $4$-cycle, and they are consistent with the conjecture that
\emph{every} trajectory with non-degenerate initial vectors, $x^0\neq y^0$ and
$(x^0,y^0)\notin E_\alpha$ converges to it. A rigorous proof of this conjecture
remains an open problem: local stability guarantees
convergence only from a neighbourhood of the cycle, and the uniqueness of
the exchange $2$-cycle (Lemma~\ref{lem_list_2cycles}) does not by itself
exclude more complicated $\omega$-limit sets, such as invariant curves or
cycles of higher periods. Partial arguments in favour of the conjecture are
the impossibility for a trajectory to remain forever in the regime of
simultaneous preservation of positions (Lemma~\ref{lem_double_active}) and
the saddle structure of the
symmetric state (Proposition~\ref{prop_saddle}), which excludes attraction
of typical trajectories to the diagonal; moreover, by
Lemma~\ref{lem_order}(1) every trajectory enters the compact set $K_\alpha$
after two steps, so that all $\omega$-limit sets lie in $K_\alpha$.
\end{remark}

\begin{figure}[htbp]
\centering
\includegraphics[width=\textwidth]{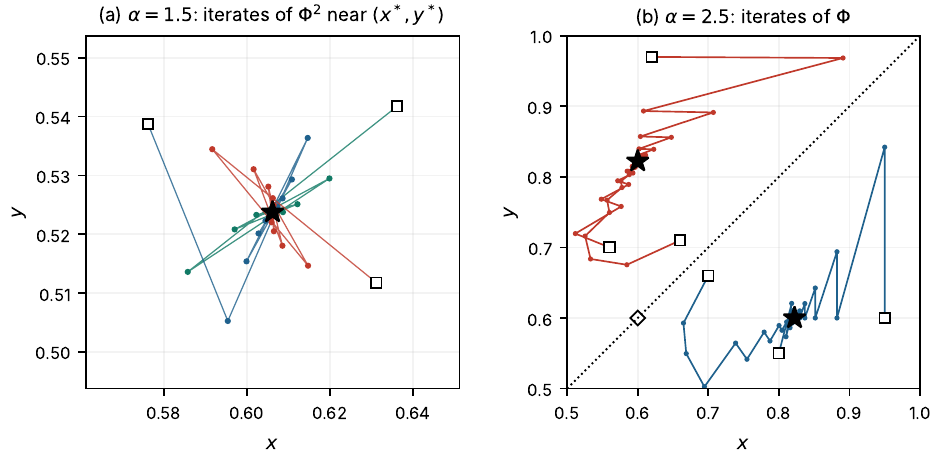}
\caption{(a) $\alpha=1.5$: iterates of $\Phi^2$ from three initial points
(squares) near the point $(x^\ast,y^\ast)$ of the exchange $2$-cycle (star).
(b) $\alpha=2.5$: trajectories of $\Phi$; those starting with $x^0>y^0$
converge to $(x^\ast,y^\ast)$, those with $x^0<y^0$ to $(y^\ast,x^\ast)$ (stars).
The diamond is the symmetric saddle $\left(1-\tfrac1\alpha,1-\tfrac1\alpha\right)$.}
\label{fig_trajectories}
\end{figure}

\begin{figure}[htbp]
\centering
\includegraphics[width=\textwidth]{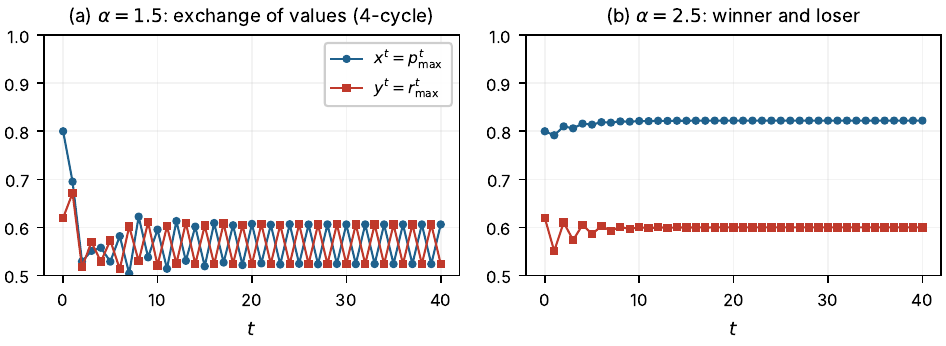}
\caption{Time series of the maxima $x^t$, $y^t$ from $(x^0,y^0)=(0.80,0.62)$.
(a) $\alpha=1.5$: the maxima exchange their values at every step (the
$4$-cycle). (b) $\alpha=2.5$: convergence to the asymmetric fixed point.}
\label{fig_time_series}
\end{figure}

Typical transients are shown in Figs.~\ref{fig_trajectories} and
\ref{fig_time_series}. Near the $4$-cycle the multipliers are complex with
argument close to $\pi$ (about $174^\circ$ for $\alpha=1.5$), so that
consecutive iterates of $\Phi^2$ lie on opposite sides of the cycle point.

\subsection{The rate of collapse at $\alpha=2$}

\begin{figure}[htbp]
\centering
\includegraphics[width=\textwidth]{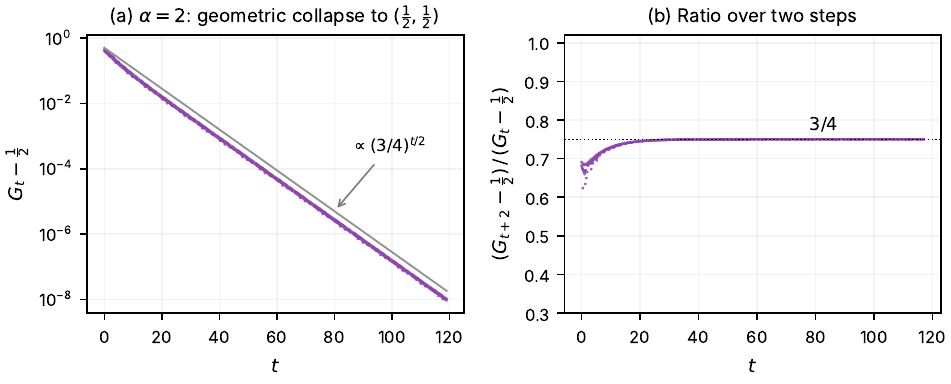}
\caption{Collapse at $\alpha=2$ for five random initial points.
(a) $G_t-\tfrac12$, where $G_t=\max(x^t,y^t)$, in logarithmic scale, with the
reference slope $\bigl(3/4\bigr)^{t/2}$ (grey). (b) The ratio
$\left(G_{t+2}-\tfrac12\right)/\allowbreak\left(G_t-\tfrac12\right)$, which tends to $\tfrac34$.}
\label{fig_collapse}
\end{figure}

Theorem~\ref{thm_alpha2} establishes convergence to $\left(\tfrac12,\tfrac12\right)$ at
$\alpha=2$ but gives no rate. Numerically, the convergence is geometric: for the five initial points of
Fig.~\ref{fig_collapse} and for $10^5$ random initial points with $x^0\neq y^0$ the
ratio $\left(G_{t+2}-\tfrac12\right)/\allowbreak\left(G_t-\tfrac12\right)$ tends to
$\tfrac34$; at $t=80$ it differs from $\tfrac34$ by less than $2\cdot10^{-6}$ for
all of them. This value coincides with the common limit of
the contraction factors of the asymmetric attractors on both sides of the
bifurcation (Fig.~\ref{fig_bif_mech}(b)).

The value $\tfrac34$ is explained by the leading-order (piecewise-linear) dynamics near the
corner; on the stability of fixed points lying on switching manifolds via
piecewise-linear approximations see \cite{simpson2020}. Write $x=\tfrac12+u$, $y=\tfrac12+v$ with $u\ge v\ge0$ small and
$\alpha=2$. Expanding \eqref{eq_macro_system}, we find that the first
coordinate is on the active branch and, up to terms of second order,
$$x'=\tfrac12+u-\tfrac v2,\qquad y'=\tfrac12+\bigl|v-\tfrac u2\bigr|,$$
the two signs of $v-\tfrac u2$ corresponding to the two branches of the
second coordinate. In terms of $s=v/u\in[0,1]$ this piecewise-linear map
reads
$$G'-\tfrac12=\left(1-\frac s2\right)\left(G-\tfrac12\right),\qquad
s'=f(s):=\frac{|s-\tfrac12|}{1-\tfrac s2},$$
and the order $x'\ge y'$ is preserved. The map $f$ sends $[0,\tfrac12]$
onto itself and is an involution there, $f(f(s))=s$; moreover, for
$s\in[0,\tfrac12]$,
$$\left(1-\frac s2\right)\left(1-\frac{f(s)}2\right)=\frac34 .$$
For $s\in\left(\tfrac12,1\right)$ one has $f(s)<s$, so that after finitely many steps
$s$ enters $[0,\tfrac12]$. Hence, to leading order, $G_t-\tfrac12$ is
multiplied by exactly $\tfrac34$ over every two steps, while the one-step
factors alternate between $1-\tfrac{s_0}{2}$ and $1-\tfrac{f(s_0)}{2}$ for
some $s_0$ depending on the trajectory (for instance, $0.864$ and $0.868$ for one of the trajectories of Fig.~\ref{fig_collapse}).
On the diagonal ($s=1$) the factor equals $\tfrac12$ at every step, in
agreement with the one-sided derivative of $g$ at $x=\tfrac12$
(Lemma~\ref{lem_diag}); thus symmetric initial states collapse faster than
asymmetric ones.
We stress that this is a leading-order computation, which agrees with the
numerics to within $2\cdot10^{-6}$ at $t=80$; a rigorous
proof would require control of the second-order terms.

\section{Conclusions}\label{concl}
We have introduced a discrete-time model of interaction concentrated at
the positions of maximal presence of two players distributed over two
positions: all fixed points and exchange $2$-cycles of its map of maxima
are found and their stability type is determined, convergence is
proved globally for symmetric initial states and for $\alpha=2$, while the global attractivity of the asymmetric attractors is
supported numerically (Section~\ref{numerics}). The
external-resource parameter $\alpha$ splits the dynamics into three regimes
separated by a single border-collision bifurcation at $\alpha=2$: an
oscillatory regime $0<\alpha<2$, with a symmetric $2$-cycle, unstable with
respect to asymmetric perturbations, and a locally stable asymmetric
$4$-cycle, governed at $\alpha=1$ by
the golden ratio and the plastic number, respectively; a global collapse to the uniform
state exactly at $\alpha=2$; and a symmetry-broken regime $\alpha>2$, in
which a winner--loser fixed point is born locally stable with multipliers
$\pm\sqrt3/2$, while the symmetric state persists only as a saddle.

From the applied standpoint, the model yields a sharp qualitative
prediction: the outcome of an interaction driven by direct competition of
leaders
is governed by a single critical level of the external resource. Below the
critical value the system sustains perpetual alternation of dominance ---
neither player can secure a lasting advantage, and the leadership positions
oscillate with an amplitude that vanishes as the external resource
approaches the critical value (Fig.~\ref{fig_bif_main}).
Exactly at the critical value the competition annihilates dominance
altogether, driving both players to the state of complete uncertainty.
Above it, an arbitrarily small initial advantage is irreversibly amplified:
the order of the maxima never changes, and the gap between them cannot
shrink to zero (Lemma~\ref{lem_order}); according to the numerical
experiments, the system then settles into a permanent winner--loser
hierarchy. The further growth of the resource accelerates the convergence
to this hierarchy and sharpens it in
relative terms: the share of the winner outside its dominant position
becomes ever smaller compared with that of the loser. The
border-collision nature of the transition means that the type of the regime
changes at the single value $\alpha=2$: the period of the attractor and the
sign of the transverse multiplier of the symmetric state change discontinuously there,
whereas the asymptotic values of the maxima pass continuously through the
uniform state. In this sense the model reproduces, in a minimal
and exactly analysable setting, both the alternation of leadership known from
technology races \cite{fudenberg1983} and the amplification of small initial
differences into dominance hierarchies \cite{bonabeau1995}. We emphasize that these conclusions concern a deliberately stylized two-position model with deterministic dynamics and a common resource parameter.

The model extends directly to a finite number $n\ge2$ of positions: at every
step the maximal coordinate of each vector is multiplied by the same factor as
in \eqref{conf_max}, while all the other coordinates change only through the
normalization. The states in which the mass of each vector outside its
dominant position is concentrated at a single position form an invariant set,
on which the dynamics coincides with the one studied here. Outside this set,
the analysis of the extended model meets three difficulties. First, the pair
of maxima no longer evolves autonomously: the new maxima depend on the
distribution of the remaining mass among the non-dominant positions through
the share of the largest non-dominant coordinate; these shares are preserved
while the dominant positions do not change and are re-evaluated at every
switching, so that the dynamics of the maxima becomes a switched system rather
than a single planar map. Second, the number of branches of the evolution law
grows with $n$, and the dominant position can move cyclically over several
positions, so that attractors without counterparts for two positions, such as
longer cycles, cannot be excluded a priori. Third, the coordinate of
${\bf p}^t$ at its dominant position increases exactly when
$\alpha(1-r^t_{\max})>1$; for two positions this requires $\alpha>2$, whereas
for $n$ positions it is possible already for $\alpha>\tfrac{n}{n-1}$, since
$r^t_{\max}\ge\tfrac1n$, so that the location of the transition can depend on
$n$. The analysis of this extension is the subject of our further research.

\section*{Data availability}
No data were used for the research described in the article.
The Python code that reproduces the numerical experiments and the figures is available at
\url{https://github.com/kseniajasko/dynamic-conflict-systems} (folder \texttt{n=2 max}).

\section*{Funding}
This work was supported by the project ``Mathematical modeling of complex systems and security-related processes'' (No.~0125U000299).

\section*{Declaration of generative AI and AI-assisted technologies in the manuscript preparation process}

During the preparation of this work the author used Claude (Anthropic) in
order to prepare and edit the English text, to check the consistency of the
notation, cross-references and statements, to propose reformulations of parts
of the text, and to write and check the code for the numerical experiments and
the figures. After using this tool, the author reviewed and edited the content
as needed and takes full responsibility for the content of the published
article.

\end{document}